\documentclass[11pt]{amsart}
\usepackage{amsmath}
\usepackage{amssymb, verbatim}
\usepackage{pstricks,pst-node,pst-plot}
\usepackage{comment}
\usepackage{color}
\usepackage{extarrows}

\usepackage[linktocpage=true,hypertexnames=false]{hyperref}
\hypersetup{colorlinks=true, citecolor=blue, linkcolor=blue, urlcolor=blue, pdfauthor={}, pdftitle={},breaklinks=true}

\usepackage[numbers]{natbib}

\title[]{The Strominger System and Flows by the Ricci Tensor}

\theoremstyle{plain}
\newtheorem{thm}{Theorem}[section]
\newtheorem{prop}[thm]{Proposition}

\newtheorem{conj}[thm]{Conjecture}

\theoremstyle{definition}
\newtheorem{ex}[thm]{Example}
\newtheorem{rk}[thm]{Remark}

\numberwithin{equation}{section}

\newcommand{\F}{\mathcal{F}}

\newcommand{\be}{\begin{equation}}
\newcommand{\bea}{\begin{eqnarray}}
\newcommand{\eea}{\end{eqnarray}} 
\newcommand{\ee}{\end{equation}}

\renewcommand{\leq}{\leqslant}
\renewcommand{\geq}{\geqslant}

\renewcommand{\epsilon}{\varepsilon}
\renewcommand{\phi}{\varphi}

\begin{document}

\author[S. Picard]{S\'ebastien Picard}
  \email{spicard@math.ubc.ca}
  \address{Department of Mathematics, The University of British
    Columbia, 1984 Mathematics Road, Vancouver BC Canada V6T 1Z2}
\thanks{The author is supported by an NSERC Discovery Grant.}

\maketitle

\begin{abstract}
This is a survey on the Strominger system and a geometric flow known as the anomaly flow. We will discuss various aspects of non-K\"ahler geometry on Calabi-Yau threefolds. Along the way, we discuss balanced metrics and balanced classes, the Aeppli cohomology class associated to a solution to the Strominger system, the equations of motion of heterotic supergravity, and a version of Ricci flow in this special geometry.
\end{abstract}


\section{Introduction}

\par K\"ahler Calabi-Yau threefolds are fundamental objects in many branches of mathematics. From the perspective of differential geometry and nonlinear PDE, a foundational result is Yau's theorem \cite{Yau78} from 1978 which equips these spaces with K\"ahler Ricci-flat metrics. In 1985, these special geometries were introduced to the field of string theory by Candelas--Horowitz--Strominger--Witten \cite{CHSW}, and have since then been central in the study of string compactifications.

An outstanding problem is to understand the parameter space of Calabi-Yau threefolds. Distinct manifolds are known to meet along interface regions by degenerations of the structure to a singular limit. There is a proposal, put forth by Reid \cite{Reid} in the algebraic geometry literature and Candelas--Green--Hubsch \cite{CGH} (see also \cite{GH1,GH2} for earlier work) in the string theory literature, which seeks to connect the space of topologically distinct threefolds into a single parameter space where two objects are connected by a sequence of degenerations and resolutions. 

It was soon noticed that there is an issue with this proposed program: the process of degeneration and resolution may connect a K\"ahler threefold to a non-K\"ahler manifold. We note here two basic examples.

\begin{ex}
This example can be found in \cite{Friedman}. Consider a smooth quintic threefold $\hat{X} \subset \mathbb{P}^4$. Take two $(-1,-1)$ curves on $\hat{X}$ and contract them: $\pi: \hat{X} \rightarrow X_0$. The resulting singular space $X_0$ has nodal singularities. Apply Friedman's theorem \cite{Friedman} (see Tian \cite{Tian} and Kawamata \cite{Kawamata} for higher order deformations) to smooth the singular points by deformation of complex structure: $X_0 \rightsquigarrow X_t$. The resulting manifold $X_t$ is complex with trivial canonical bundle, but it cannot admit a K\"ahler structure. This is because $b_2(\hat{X})=1$, so upon contraction of a generator of homology the resulting manifold has $b_2(X_t)=0$. The process $\hat{X} \rightarrow X_0 \rightsquigarrow X_t$ of holomorphic contraction followed by deformation of complex structure is sometimes called a conifold transition, and we refer to Rossi's survey \cite{Rossi} for a detailed introduction.
\end{ex}

\begin{ex}
The previous example starts from a quintic and degenerates holomorphic curves. We may alternatively degenerate the complex structure of a smooth quintic by varying the coefficients of the polynomial cutting out the space in $\mathbb{P}^4$. It is possible to tune coefficients such that a smooth quintic degenerates to a nodal quintic with a single node (see e.g. Chapter D of \cite{Hubsch}). A small resolution of this singular space cannot be K\"ahler (see e.g. Chapter 4 of \cite{Hubsch} or Section 6 of \cite{COGP}). This process $X_t \rightsquigarrow X_0 \leftarrow \hat{X}$ is a deformation of complex structure followed by small resolution, and this sort of operation is sometimes called a reverse conifold transition.
\end{ex}

Since degenerations and resolutions do not preserve the K\"ahler condition, the idea of Friedman \cite{Friedman} and Reid \cite{Reid} is to include the limiting spaces into the parameter space of threefolds. Said otherwise, we expand our definition of the Calabi-Yau threefold to include those non-K\"ahler complex manifolds which can be reached by degenerations. These spaces still inherit some of the structure of K\"ahler Calabi-Yau threefolds; for examples the stability of the tangent bundle is preserved through conifold transitions \cite{CPY,FLY}, and for a discussion of the geometry of the reverse process see \cite{STY,GiuSpo}.

From the perspective of differential geometry, the question, emphasized by S.-T. Yau throughout the years (e.g. \cite{LY05,YauNadis}), is how to geometrize the objects connected to K\"ahler threefolds by degenerations such as a conifold transition. The idea is that there should be a uniform geometric structure which holds whether the underlying space admits a K\"ahler structure or not. The first milestone along this direction is the work of Fu--Li--Yau \cite{FLY}, which solves the non-K\"ahler constraint $d \omega^2=0$ through a conifold transition.

The proposed set of equations to geometrize conifold transitions are as follows. Inspiration comes from supersymmetric string theory \cite{Strominger}. For a compact complex manifold $X$ of dimension 3 with holomorphic volume form $\Omega$, the equations of heterotic string theory for the hermitian metric $\omega$ are:
  \[
d(|\Omega|_\omega \, \omega^2)=0, \quad i \partial \bar{\partial} \omega = \alpha' ({\rm Tr} \, R \wedge R - {\rm Tr} \, F \wedge F).
  \]
  Here $\alpha'>0$, and we note that it can be set to any positive constant by rescaling $\omega \mapsto \lambda \omega$. The form $F \in \Lambda^{1,1}({\rm End} \, E)$ is the curvature form of a connection on a vector bundle $E \rightarrow X$ which is Hermitian-Yang-Mills with respect to $\omega$. K\"ahler Calabi-Yau geometries are special solutions to the system with $E=T^{1,0}X$ and $F=R$. We note that Fu--Li--Yau \cite{FLY} solve the first equation through conifold transitions, and the appearance of a bundle $E$ in the second equation through Calabi-Yau degenerations such as conifold transitions is not yet understood, though we refer to e.g. \cite{ABG, CPY, Chuan} for recent developments.

Yau's proposal on the geometrization of conifold transitions \cite{LY05,FLY,CPY,YauNadis} is to equip all threefolds connected by degenerations and resolution with the non-K\"ahler structure given by Strominger's equations for suitable bundles $E$. The geometrization of conifold transitions remains an open problem which will not be further elaborated on in this survey. Instead, we take this background discussion as context and motivation for studying the differential geometry and PDE aspects of the non-K\"ahler geometry of Strominger's equations on a complex manifold.

  We will start by presenting the basic features of the equations and discuss how to frame the problem as a search for special representatives in a chosen cohomology class. The second part of this article will survey recent works on a parabolic version of the Strominger system. This geometric flow, introduced in \cite{PPZ18} as the anomaly flow, is a version of Hamilton's Ricci flow \cite{Hamilton} with additional non-K\"ahler terms.
  \smallskip
  \par {\bf Acknowledgements:} The author thanks D.H. Phong and S.-T. Yau for countless discussions on the Strominger system throughout the years. The author also thanks T.C. Collins, T. Fei, J. McOrist, E. Svanes, P.-L. Wu and X.-W. Zhang for collaborations on the works presented in this survey, and B. Friedman, M. Garcia-Fernandez, and C. Suan for helpful discussions. The author thanks T. Hubsch, C. Shabazi and S. Sethi for comments.

  \section{Overview of the system}
\subsection{Introducing the equations}
 We start by establishing notation. Let $X$ be a compact complex manifold of complex dimension 3. Let $\Omega$ be a nowhere vanishing  holomorphic $(3,0)$ form. Let $E \rightarrow X$ be a complex vector bundle. A hermitian metric on $X$ will be denoted $\omega = i g_{\mu  \bar{\nu}} \, dz^\mu \wedge d z^{\bar{\nu}}$. A metric on the bundle $E$ will be denoted $h$. The Strominger system \cite{Strominger} is for an unknown pair $(\omega,h)$ solving
\begin{align} \label{confbal}
  & d (|\Omega|_\omega \, \omega^2)=0, \\
  & F_h \wedge \omega^2 =0, \label{hym} \\  
  & i \partial \bar{\partial} \omega = \alpha' ({\rm Tr} \, R_\omega \wedge R_\omega - {\rm Tr} \, F_h \wedge F_h ). \label{ac}
\end{align}
Here $F_h$ and $R_\omega$ are the Chern connections of $h$ and $g$. 

The first equation \eqref{confbal} is the conformally balanced equation. It is, after a conformal change, the balanced condition in non-K\"ahler complex geometry studied by Michelsohn \cite{Michelsohn}, which is
\[
d \omega^{n-1} = 0
\]
on a complex manifold of dimension $n$. This condition is dual to the K\"ahler condition, as the K\"ahler condition is $d \omega = 0$ while the balanced condition is $d \star \omega = 0$. Not every complex manifold admits a balanced metric, though there are many non-K\"ahler examples \cite{Michelsohn}. It is required in heterotic string theory on $\mathbb{R}^{3,1} \times X$ as the conformally balanced condition is related to a parallel spinor generating a supersymmetry \cite{Strominger} (see \cite{MGFRT-ANN} for an alternate proof).

The second equation \eqref{hym} is the Hermitian-Yang-Mills equation
\[
F_h \wedge \omega^{n-1} = \mu \, {\rm id}_E \, \omega^n
\]
in the non-K\"ahler setting. Its solvability as a stand-alone equation follows from the the Li-Yau \cite{LY86} extension of the Donaldson-Uhlenbeck-Yau theorem \cite{Donaldson,UY}.  This theorem states that Hermitian-Yang-Mills metrics exist if and only if the bundle $E \rightarrow (X,\omega)$ is polystable. In the setup discussed here, we assume $c_1(E)=0$ so that the constant $\mu$ is zero, as can be seen by taking the trace and integrating.

We will refer to \eqref{ac} as the anomaly cancellation equation, though in string theory it is also called the 4-form heterotic Bianchi identity. It implies the Chern class condition $c_2^{\rm BC}(X)=c_2^{\rm BC}(E)$ in Bott-Chern cohomology. A discussion and definition of Bott-Chern cohomology will be given in \S \ref{section:balanced}. As this equation is perhaps the least familiar in the mathematics literature, we give an outline of the interpretation of this equation as a sort of holomorphic structure on
\[
Q = T^{*(1,0)}X \oplus {\rm End} \, E \oplus T^{1,0}X.
\]
The construction we give here is joint work with J. McOrist and E. Svanes \cite{MPS}, though this idea without the ${\rm Tr} \, R_\omega \wedge R_\omega$ term with Chern connection was introduced by de la Ossa-Svanes \cite{dlOSvanes} and related versions of this construction in the mathematics literature can be found in works of Bismut \cite{Bismut} and Gualtieri \cite{Gualti}. Consider the operator 
\[
  \bar{D}: \Lambda^{0,k}(Q) \rightarrow \Lambda^{0,k+1}(Q)
\]
given by
\[
  \bar{D} = \begin{bmatrix} \bar{\partial} & -2 \alpha' \F^* & \mathcal{H} - 2 \alpha' \mathcal{R} \nabla \\
    0 & \bar{\partial}_A & \F \\
  0 & 0 & \bar{\partial} \end{bmatrix},
\]
acting on
\[
\begin{bmatrix} \eta \\ \phi \\ \mu \end{bmatrix} \in  \Gamma \begin{bmatrix} T^{*(1,0)}X \\ {\rm End} \,E \\ T^{1,0}X \end{bmatrix} \otimes \Lambda^{0,k}.
\]
The definitions are
\begin{align*}
  & \F: \Lambda^{0,k}(T^{1,0}X) \rightarrow \Lambda^{0,k+1}({\rm End} \, E), \quad \F(\mu) = F_{\alpha \bar{\beta}} \, dz^{\bar{\beta}} \wedge \mu^\alpha, \nonumber\\
  & \F^*: \Lambda^{0,k}({\rm End} \, E) \rightarrow \Lambda^{0,k+1}(T^{*(1,0)}), \quad \F^*(\phi) = {\rm Tr} \, F_{\alpha \bar{\beta}} \, dz^\alpha \otimes d z^{\bar{\beta}} \wedge \phi \nonumber\\
  & \mathcal{H} : \Lambda^{0,k}(T^{1,0}X) \rightarrow \Lambda^{0,k+1}(T^{*(1,0)}), \quad \mathcal{H}(\mu) = H_{\alpha \bar{\beta} \lambda} \, dz^\alpha \otimes dz^{\bar{\beta}} \wedge \mu^\lambda.
\end{align*}
Here $H = i (\partial-\bar{\partial})\omega$. The operator
\[
  \mathcal{R} \nabla:  \Lambda^{0,k}(T^{1,0}X) \rightarrow \Lambda^{0,k+1}(T^{*(1,0)}),
\]
involves the Strominger-Bismut \eqref{Str-Bis} connection $\hat{\nabla}$ and is
\[
  (\mathcal{R} \nabla) (\mu) = - {1 \over k!} R_{\rho \bar{\beta}}{}^\sigma{}_\lambda \hat{\nabla}_\sigma \mu^\lambda{}_{\bar{\kappa}_1 \dots \bar{\kappa}_k} dz^\rho \otimes dz^{\bar{\beta} \bar{\kappa}_1 \dots \bar{\kappa}_k}.
\]
The calculation of \cite{MPS} establishes the identity
\[
\bar{D}^2 \begin{bmatrix} \eta \\ \phi \\ \mu \end{bmatrix}  = \begin{bmatrix} {1 \over 2}( - i \partial \bar{\partial} \omega - \alpha' {\rm Tr} \, F \wedge F + \alpha' {\rm Tr} \, R \wedge R)_{\alpha \bar{\beta} \gamma \bar{\delta}} dz^\alpha \otimes dz^{\bar{\beta} \bar{\delta}} \wedge \mu^\gamma \\ 0 \\ 0 \end{bmatrix}.
\]
Thus the nilpotency condition $\bar{D}^2 = 0$ is equivalent to the anomaly cancellation equation \eqref{ac}. This suggests that we could interpret the anomaly cancellation equation as a holomorphic structure on the bundle
\[
(Q, \bar{D}) \rightarrow X.
\]
However, $\bar{D}$ is not a traditional holomorphic structure due to the $2 \alpha' \mathcal{R} \nabla$ term. It is a 1st order differential operator, which is a holomorphic $\bar{\partial}$ operator plus an off-diagonal $\alpha'$-correction.

\begin{rk}
  We note that this idea of relating the anomaly cancellation equation to a holomorphic structure on $Q$ can be taken much further when ${\rm Tr} \, R \wedge R$ is taken with respect to an auxiliary Hermitian-Yang-Mills connection $\nabla$ on $T^{1,0}X$ (see discussion in \S \ref{section-alt} below). In that case, there is an associated operator $\bar{D}$ which is a true holomorphic structure. In \cite{MGFRT-ANN,MGFRT-TAMS} the holomorphic bundle $Q$ is also equipped with a Dorfman bracket structure
  \[
[ \cdot, \cdot ]: Q \otimes Q \rightarrow Q
  \]
  making $Q$ a holomorphic Courant algebroid. Thus the anomaly cancellation equation can be interpreted as a type of Courant algebroid structure \cite{MGFRT-TAMS}.

  This link to generalized geometry has led to various insights \cite{ASTW,ADMSE,MGFRST,MGFRT-ANN,MGFRT-JDG} into the geometry of the system and its moduli. Including both the anomaly cancellation equation and the conformally balanced equation leads to a Hermitian-Yang-Mills equation on $Q$ \cite{dlOLS}, and from there Garcia-Fernandez and Gonzalez Molina \cite{MGFM23a,MGFM23b} define a Futaki invariant and conjecture obstructions and stability notions for the solvability of the system. We refer to \cite{MPS} for the Hermitian-Yang-Mills interpretation of the system with Chern connection. From the perspective of quantum physics, these ideas of structures on $Q$ lead to a calculation of the heterotic 1-loop partition function \cite{AIMSSTW}.
\end{rk}

\subsection{Examples} The Strominger system \eqref{confbal}, \eqref{hym}, \eqref{ac} with $\alpha'=0$ requires $\omega$ to be K\"ahler Calabi-Yau (e.g. \cite{Ivanovpapa}) and $E \rightarrow (X,\omega_{\rm CY})$ a stable holomorphic bundle. When $\alpha'>0$, K\"ahler Calabi-Yau metrics still provide solutions to the system by letting $E= T^{1,0}X$ and $h = g_{\rm CY}$. We briefly note here in passing the three types of known compact non-K\"ahler solutions in the literature:

\begin{ex}
  Perturbations from K\"ahler Calabi-Yau threefolds. This approach was initiated by Li--Yau \cite{LY05}, and further developed by Andreas and Garcia-Fernandez \cite{AMGF}. In this setup, $E \rightarrow (X,\omega_{\rm CY})$ is a stable degree zero holomorphic vector bundle over a K\"ahler Calabi-Yau manifold $(X,\omega_{\rm CY})$ with $c_2(E)=c_2(X)$. Solutions for small $\alpha'>0$ are obtained by analysis of the linearized operator and implicit function theorem. An alternate proof is also presented in \cite{CPY-IFT}.
\end{ex}

\begin{ex}
  Non-K\"ahler fibrations. In this approach, $X$ is the total space of a fibration and the equations are reduced to the base of the fibration. There are two known reductions of this type.

  \begin{itemize}
  \item The first type of fibration $\pi: X \rightarrow K3$ is a $T^2$ fibration over a $K3$ surface. This sort of construction of a non-K\"ahler manifold goes back to Calabi--Eckmann \cite{CE} and was suggested for study in string theory in \cite{DRS}. Fu--Yau \cite{FY} solved the Strominger system \eqref{confbal}, \eqref{hym}, \eqref{ac} by reducing the equations to a Monge-Amp\`ere equation on the base $K3$. For related work on these non-K\"ahler threefolds, we refer to e.g. \cite{BecSeth, BBFTY, FGV, MGF-T, GP, MMS}.
  \item The second type of fibration construction is developed in \cite{Fei,FHP}, where solutions are found on $T^4$ fibrations $\pi: X \rightarrow \Sigma$ over a genus $g \geq 2$ Riemann surface $\Sigma$. The Strominger system is solved by reducing the equations to a nonlinear PDE on a Riemann surface.
    \end{itemize}
\end{ex}

\begin{ex}
Lie group constructions. Here $X$ is assumed to be a homogeneous space with action by a Lie group. For solutions with symmetry, see e.g. \cite{FeiYau,FIUV,Grant,OUV}.
\end{ex}

Thus there are many known examples of solutions to the Strominger system, and this is due to the fact that the non-K\"ahler structure allows for more holomorphic symmetries than allowed for torsion-free Calabi-Yau structures. The class of manifolds solving the system is strictly greater than K\"ahler Calabi-Yau threefolds and in fact, the solutions in \cite{FHP} admit infinitely many topological types.

\subsection{Alternate formulations} \label{section-alt}
Before moving on, we note here that there are three common variants of the system of equations presented in this survey:
\smallskip
\par $\bullet$ System with Chern connection. These are the equations \eqref{confbal}, \eqref{hym}, \eqref{ac} where the unknowns are $(\omega,h)$ and the relevant curvature forms are
\[
F_h=\bar{\partial}(h^{-1} \partial h), \quad R_\omega = \bar{\partial}(g^{-1} \partial g).
\]
In other words, $R$ and $F$ are the curvature forms of the Chern connections of $g$ and $h$. This is the system solved and studied in e.g. \cite{Fei,FHP,FTY,LY05,FY,CPY-IFT,PW,PPZ18}. We refer to this system as the Strominger system \cite{Strominger}.
\smallskip
\par $\bullet$ System with auxiliary Hermitian-Yang-Mills connection. These equations are for an unknown triple $(\omega,\nabla,A)$ satisfying:
\begin{align*}
& i \partial \bar{\partial} \omega = \alpha' ({\rm Tr} \, R_\nabla \wedge R_\nabla - {\rm Tr} \, F_A \wedge F_A),\quad d (|\Omega|_\omega \, \omega^2)=0, \\
& F_A^{0,2}=0, \quad F_A \wedge \omega^2 =0, \quad R_\nabla^{0,2}=0, \quad R_\nabla \wedge \omega^2 = 0.
\end{align*}
Here $R$, $F$ are the curvature forms of unitary connections $A$,$\nabla$ on the bundles $E$, $T^{1,0}X$. This system is studied and solved for various connections $\nabla$ in e.g. \cite{AMGF,ADMSE,dlOSvanes,Ivanov,MGF-Survey,MGF-T,MGFM23b}. Compared to the system with Chern connection, there is an additional unknown Hermitian-Yang-Mills connection $\nabla$ on $T^{1,0}X$, which has the effect of introducing a symmetry between $F_A$ and $R_\nabla$.
\smallskip
\par $\bullet$ System with Hull connection. The statement of constraints is
\begin{align*}
  &  i \partial \bar{\partial} \omega = \alpha' ({\rm Tr} \, R^{\rm H} \wedge R^{\rm H} - {\rm Tr} \, F \wedge F) + O(\alpha'^2), \\
&  F^{0,2}=O(\alpha'), \quad F \wedge \omega^2 =O(\alpha'), \quad  d (|\Omega|_\omega \, \omega^2)=O(\alpha'^2).
\end{align*}
Here $R^{\rm H}$ is taken with respect to the Hull connection, whose definition we recall in Appendix \ref{appendix}. As $R^{\rm H}$ is not a $(1,1)$-form, we cannot expect equality to hold in the $i \partial \bar{\partial} \omega$ equation without additional $O(\alpha'^2)$ contributions since ${\rm Tr} \, R^{\rm H} \wedge R^{\rm H}$ will not be a $(2,2)$-form. In string theory these equations come from an expansion in a small parameter $\alpha'$, and the equations are derived order by order in the expansion. The string theoretic reason for selecting the Hull connection is explained in \cite{Hull} and more recently in \cite{MMS}. Works on this system include e.g. \cite{AMP, CMOS, CCDL,Hull, MartelliSparks,McOSva,MMS}. 

Each formulation comes with various distinct features of interest. The system with Chern connection has fewer equations and unknowns and the nonlinearity in ${\rm Tr} \, R \wedge R$ makes its analysis interesting as a fully nonlinear geometric PDE. The system with auxiliary connection allows a geometric reinterpretation in terms of structures in generalized geometry \cite{ASTW,MGFRT-ANN,MGFRT-JDG}. The heterotic system with Hull connection is consistent with various perspectives in string theory and matches calculations from both the spacetime supersymmetry and the worldsheet supersymmetry points of view \cite{Hull, MMS}.

In summary, these different interpretations of the system have attracted various mathematical communities including researchers in the fields of analysis and PDE \cite{FY, FLY, LY05, PPZ18CAG}, generalized geometry \cite{MGF-Survey, MGFRST, MGFRT-JDG}, and string theory \cite{AGS, ADMSE, BecSeth, BBDG, BBDGS, CMO, dlOSvanes, GMW}.

In the current survey, we will focus our attention on the Strominger system with Chern connection to streamline presentation, but we will also sometimes mention related results on the other formulations along the way.


\section{Cohomology classes} \label{section:coho}

\subsection{Overview}
There is a theme in differential geometry to find optimal representatives in a given cohomology class. This section will explore analogs of the K\"ahler class for non-K\"ahler Calabi-Yau threefolds and their special representatives. For context, we recall two fundamental theorems on special representatives and cohomology.

\begin{ex}
In Hodge theory, we start from a de Rham cohomology class $[\alpha] \in H^k(X,\mathbb{R})$. The Hodge theorem states that on a compact Riemannian manifold $(M,g)$, every cohomology class contains a unique harmonic representative
\[
\alpha_{\rm H} \in [\alpha]
\]
solving $\Delta_g \alpha_{\rm H} = 0$, where $\Delta_g = d d^\dagger + d^\dagger d$ and $d^\dagger$ is the $L^2$ adjoint associated to the metric $g$.
\end{ex}

\begin{ex}
  Let $(X,\omega)$ be a compact K\"ahler Calabi-Yau manifold. From a K\"ahler metric $\omega$, since
  \[
d \omega = 0
  \]
 one associates the K\"ahler class
\[
[\omega] \in H^{1,1}(X,\mathbb{R}).
\]
Yau's theorem \cite{Yau78} states that every K\"ahler class contains a unique representative
\[
\omega_{\rm CY} \in [\omega]
\]
which is K\"ahler Ricci-flat. Thus on a K\"ahler Calabi-Yau manifold with fixed complex structure, families of K\"ahler Ricci-flat metrics on $X$ are parametrized by a path of K\"ahler classes in $H^{1,1}(X,\mathbb{R})$.

\end{ex}

Returning to non-K\"ahler Calabi-Yau threefolds, we seek an analog of the K\"ahler class. The hermitian metric $\omega$ is no longer closed $d \omega \neq 0$, so it does not define an element in de Rham cohomology. However, there are other cohomology groups on a non-K\"ahler complex manifold involving the $\partial$ and $\bar{\partial}$ operators: Dolbeault cohomology, Bott-Chern cohomology, and Aeppli cohomology (see e.g. \cite{PopoviciBook} for an introduction). From a solution to the Strominger system one can associate two classes in non-K\"ahler cohomology:

\begin{itemize}
\item The first is the balanced class.
\[
\mathfrak{b} = [|\Omega|_\omega \, \omega^2] \in H^{2,2}_{\rm BC}(X,\mathbb{R}).
\]
One could setup the problem of solving the Strominger system as looking for solutions in a given Bott-Chern balanced class.

  \item The second is the Aeppli class
\[
\mathfrak{a} = [\omega - \alpha' R_2[g,\hat{g}]+ \alpha' R_2[h,\hat{h}] - \alpha' \hat{\beta}] \in H^{1,1}_{\rm A}(X,\mathbb{R}).
\]
Alternatively, one could look for solutions to the Strominger system in a given Aeppli class.
  \end{itemize}

  The precise definitions of the above classes will be given in \S \ref{section:balanced} and \S \ref{section:aeppli} below. The broad goal will be to formulate existence and uniqueness theorems in a given cohomology class, and use these classes as parameters for families of solutions to the Strominger system. 

  Before expanding on the details, we note that on a K\"ahler Calabi-Yau threefold there are isomorphisms
\[
  H^{2,2}_{\rm BC}(X,\mathbb{C}) \cong H^{2,2}_{\bar{\partial}}(X,\mathbb{C}) \cong H^{1,1}_{\bar{\partial}}(X,\mathbb{C}), \quad   H^{1,1}_{\rm A}(X,\mathbb{C}) \cong H^{1,1}_{\bar{\partial}}(X,\mathbb{C}).
\]
Thus the balanced and Aeppli class of a solution to the Strominger system can be understood as non-K\"ahler analogs of the K\"ahler class of a Calabi-Yau metric. We note that the single notion of a K\"ahler class breaks into two when extending the concept from the K\"ahler to the non-K\"ahler setting, though there is a duality in cohomology via the pairing
\[
H^{1,1}_{\rm A}(X,\mathbb{C}) \times H^{2,2}_{\rm BC}(X,\mathbb{C}) \rightarrow \mathbb{C}, \quad (\mathfrak{a}, \mathfrak{b}) \mapsto \int_X \mathfrak{a} \wedge \mathfrak{b}.
\]
There is thus the possibility of a duality of solutions to the Strominger system in a Bott-Chern balanced class or Aeppli class. At present, the Bott-Chern balanced class and Aeppli class approaches to the Strominger system detailed below are not evidently related.
 
\subsection{Balanced class} \label{section:balanced}
Let $X$ be a complex manifold of dimension $3$ with nowhere vanishing holomorphic 3-form $\Omega \in \Lambda^{3,0}(X,\mathbb{C})$.
\[
d \Omega = 0 .
\]
Let $\omega$ be a hermitian metric on $X$ satisfying the closedness condition
\begin{equation} \label{def-conf-bal}
d (|\Omega|_\omega \, \omega^2) = 0.
\end{equation}
Our conventions are $\omega= i g_{\mu \bar{\nu}} \, dz^\mu \wedge d z^{\bar{\nu}}$ and
\begin{equation} \label{def-norm-Omega}
|\Omega|_\omega^2 = {f \bar{f} \over \det g_{\mu \bar{\nu}}}, \quad \Omega = f(z) \, dz^1 \wedge dz^2 \wedge d z^3.
\end{equation}
Metrics solving \eqref{def-conf-bal} will be called conformally balanced. This structure defines a cohomology class
\[
\mathfrak{b}(\omega) = [|\Omega|_\omega \, \omega^2] \in H^{2,2}_{\rm BC}(X,\mathbb{R}),
\]
where
\[
H^{2,2}_{\rm BC}(X,\mathbb{R}) = \frac{(\ker d) \cap \Lambda^{2,2}(X,\mathbb{R})}{{\rm Im} \, i \partial \bar{\partial}}.
\]
The equation \eqref{def-conf-bal} also defines a de Rham cohomology class in $H^4_{\rm dR}(X)$, but in non-K\"ahler complex geometry the Bott-Chern cohomology $H^{2,2}_{\rm BC}(X)$ encodes more refined information. If $X$ satisfies the $\partial \bar{\partial}$-lemma, then de Rham and Bott-Chern cohomology coincide.

\subsubsection{Balanced classes and the Gauduchon conjecture}
After a conformal change $\tilde{\omega} = |\Omega|_\omega^{1/2} \omega$, the equation \eqref{def-conf-bal} becomes the balanced condition $d \tilde{\omega}^{n-1}=0$ in dimension $n=3$. The systematic study of balanced metrics was initiated by Michelsohn \cite{Michelsohn}, and the condition appeared early in complex geometry (e.g. \cite{Gauduchon75}) as a dual condition to the K\"ahler property. 

We move on to the problem of finding special representatives in a given balanced class. Namely, we seek an analog of Yau's theorem \cite{Yau78} on non-K\"ahler Calabi-Yau threefolds. There are two main conjectures in the literature along these lines. The first is the balanced version of the Gauduchon conjecture.

\begin{conj} \cite{FWW,Tos15}
Let $X$ be a compact complex manifold of dimension $n$ with holomorphic volume form $\Omega$. Let $\omega$ be a hermitian metric satisfying the balanced condition $d \omega^{n-1}=0$ and denote its balanced class by $ [\omega^{n-1}] \in H_{\rm BC}^{n-1,n-1}(X,\mathbb{R})$. There exists a balanced metric $\tilde{\omega}$ with $\tilde{\omega}^{n-1} \in [\omega^{n-1}]$ such that ${\rm Ric}_{\tilde{\omega}}=0$. 
  \end{conj}

This conjecture is known when the background manifold $X$ admits a K\"ahler metric \cite{TW17} and there other constructions in \cite{FWW,FWW2,GiuSpo}, but otherwise this is an open problem in general. Note that even if $X$ admits a background K\"ahler metric, a given balanced class $[\omega^{n-1}]$ does not necessarily contain a K\"ahler metric \cite{FuXiao}. We should clarify that here ${\rm Ric}_{\tilde{\omega}}$ denotes
  \[
i {\rm Ric}_{\tilde{\omega}} = i \partial \bar{\partial} \log |\Omega|^2_{\tilde{\omega}},
  \]
  and for non-K\"ahler metrics $\omega$ this differs from the Riemannian Ricci curvature of the metric $g$.

\subsubsection{Balanced classes and the Strominger system}
  We return now to the Strominger system, which suggests an alternate approach to special representatives in a given balanced class inspired by string theory. From the point of view of string theory, conformally balanced metrics already all satisfy a non-K\"ahler analog of Ricci flatness. Indeed, the following result is shown in \cite{FinoGrant}: let $(X,\Omega)$ be a compact complex manifold with holomorphic volume form. Then $\omega$ satisfies the conformally balanced constraint \eqref{def-conf-bal} if and only if
\begin{equation} \label{Bismut-Ricci-flat}
\hat{R}_{mn}{}^\mu{}_\mu = 0.
\end{equation}
Here are the relevant definitions. The curvature tensor appearing in \eqref{Bismut-Ricci-flat} is the cuvature of the Strominger-Bismut connection $\hat{\nabla}$ \eqref{Str-Bis}, and our cuvature conventions are detailed in \eqref{curvature}. We recall that $m,n$ denote real indices while $\mu$ denotes holomorphic indices. Thus the conformally balanced equation has the geometric interpretation of vanishing Strominger-Bismut Ricci curvature.

There is another interpretation in terms of holonomy groups. On a compact complex manifold with trivial canonical bundle, the conformally balanced equation is equivalent to the holonomy constraint
\[
{\rm Hol}(\hat{\nabla}) \subseteq SU(n).
\]
We refer to e.g. \cite{MGF-Survey} for a proof of this fact. Compared to the standard setup for special holonomy problems in differential geometry, the Levi-Civita connection is replaced by the Strominger-Bismut connection $\hat{\nabla}$, and the statement is an analog of the fact that K\"ahler Calabi-Yau manifolds have Riemannian holonomy contained in $SU(n)$.

From the point of view of the Strominger system, a special representative in a given balanced class would solve the remaining equations in the system. Concretely, we would substitute the ansatz
\begin{equation} \label{balanced-ansatz}
|\Omega|_{\omega_\Theta} \, \omega_{\Theta}^2 = |\Omega|_\omega \, \omega^2 + \Theta > 0, \quad \Theta \in \Lambda^{2,2} \cap {\rm Im} \, \partial \bar{\partial}
\end{equation}
into the remaining equations in the system,
\begin{equation} \label{balanced-remaining}
i \partial \bar{\partial} \omega_\Theta = \alpha' ({\rm Tr} \, R_\Theta \wedge R_\Theta - {\rm Tr} \, F_h \wedge F_h), \quad F_h \wedge \omega_\Theta^2=0,
\end{equation}
and solve for the unknown $(2,2)$-form $\Theta$ and bundle metric $h$. Here $R_\Theta$ denotes the Chern curvature of the metric $\omega_\Theta$. This setup would produce a solution to the Strominger system inside a given balanced class.
\[
|\Omega|_{\omega_\Theta} \, \omega_\Theta^2 \in [|\Omega|_\omega \,\omega^2].
\]
We should explain why this is an elliptic equation in the unknown $\Theta$. First, the ansatz \eqref{balanced-ansatz} does define a metric $\omega_\Theta>0$ by the square root construction for positive $(n-1,n-1)$ forms (see Michelsohn \cite{Michelsohn} or \cite{PPZ18} in the conformally balanced case). We will not use the explicit formula for $\omega_\Theta$, and instead use the simpler formula for its variation: if
\[
|\Omega|_\omega \, \omega^{2} = \Psi,
\]
then the variation of this equation works out to be \cite{PPZ18CAG} (see also \cite{PW} for an alternate presentation)
\begin{equation} \label{variation-omega}
\delta \omega = {1 \over 2 |\Omega|_\omega} \Lambda_\omega \delta \Psi.
\end{equation}
To linearize \eqref{balanced-remaining} in the direction $\delta \Theta$, we hold $h$ constant and vary $\Theta$. Substituting \eqref{variation-omega} gives
\[
i \partial \bar{\partial} \bigg[ {1 \over 2 |\Omega|_\omega} \Lambda_\omega \delta \Theta \bigg] = 2 \alpha' {\rm Tr} \, R_\Theta \wedge \delta R_\Theta.
\]
It is well-known that
\begin{equation} \label{derivative-R}
\delta R_g = \bar{\partial} \partial_\nabla (g^{-1} \delta g).
\end{equation}
Next, for a $(2,2)$-form $\eta$, we have the non-K\"ahler identity (e.g Chapter VI Theorem 6.8 in \cite{BigDemailly}) given by
\[
[i \Lambda, \bar{\partial}] \eta = - \partial^\dagger \eta + H * \eta
\]
where $H*\eta$ denotes various contractions of the tensor $H= i (\partial-\bar{\partial})\omega$ with $\eta$. Therefore the linearized equation of \eqref{balanced-remaining} is, up to lower order terms, given by
\[
{1 \over 2 |\Omega|_\omega} \Delta_\partial  (\delta \Theta) - {\alpha' \over |\Omega|_\omega}  {\rm Tr} \, R_\Theta \bar{\partial} \partial_\nabla (g_\Theta^{-1} \Lambda \delta \Theta) + l.o.t. = 0
\]
where $\Delta_\partial = \partial \partial^\dagger + \partial^\dagger \partial$, and we used that $\Theta$ and $\delta \Theta$ are closed. The ellipticity of $\Delta_\partial$ may be used to absorb the $\alpha'$-corrections when $\alpha' R_{\Theta}$ is small. Thus there exists an $\epsilon>0$ such that the equation \eqref{balanced-remaining} is elliptic in the regime
\[
|\Omega|_\omega \, \omega^2 + \Theta > 0, \quad |\alpha' R_\Theta| < \epsilon.
\]
As \eqref{balanced-remaining} is a fully-nonlinear PDE for the unknown $\Theta$, its ellipticity region is an open set rather than the entire space of $(2,2)$ forms.

To obtain the full linearized operator, we should also vary $h$ and represent the linearized operator of the system \eqref{balanced-remaining} as a $2 \times 2$ matrix acting on variations $(\delta \Theta, \delta h)$. The property that the linearized operator is elliptic allows for regularity results and perturbation theorems, and we refer to \cite{CPY-IFT, PW} for applications. Returning to the overall vision for this approach to the Strominger system via balanced classes, the main conjecture for the existence of solutions is as follows:

\begin{conj} \label{yauconj} (Yau's conjecture for balanced classes \cite{Yauconj})
  Let $X$ be a compact complex manifold of dimension $3$ with holomorphic volume form $\Omega$. Let $\omega$ be a hermitian metric satisfying the conformally balanced condition $d ( | \Omega |_\omega \, \omega^2)=0$ and denote its balanced class by $\mathfrak{b}= [|\Omega|_\omega \, \omega^{2}] \in H^{2,2}(X,\mathbb{R})$. Let $E \rightarrow X$ be a stable holomorphic bundle with respect to $\mathfrak{b}$ such that $c_1^{\rm BC} (E)=0$, $c_2^{\rm BC}(E)=c_2^{\rm BC}(X)=0$. Then there exists a parameter $\alpha'>0$ and a pair of metrics $(\tilde{\omega},h)$ such that $|\Omega|_{\tilde{\omega}} \, \tilde{\omega}^2 \in \mathfrak{b}$ and
  \[
  i  \partial \bar{\partial} \tilde{\omega} = \alpha' ( {\rm Tr} \, R_{\tilde{\omega}} \wedge R_{\tilde{\omega}} - {\rm Tr} \, F_h \wedge F_h ), \quad F_h \wedge \tilde{\omega}^2=0.
  \]
  \end{conj}

  We note here that the conjecture stated above is slightly stronger than the one stated in \cite{Yauconj} as the balanced class of the solution is prescribed. The statement with balanced class specified is called the strong version of Yau's conjecture in \cite{MGFM23a}. However, we note that our setup here differs slightly from \cite{MGFM23a}, as we require ${\rm Tr} \, R \wedge R$ to be computed with respect to the Chern connection of $\omega$, rather than an auxiliary Hermitian-Yang-Mills connection $\nabla$.
 
The conjecture is known to be true when the balanced class comes from a K\"ahler class $[\omega_{\rm CY}]$. This is joint work with T. Collins and S.-T. Yau \cite{CPY-IFT}. There exists $\epsilon>0$ such that for all $0< \alpha'<\epsilon$, we find a solution $(\omega,h)$ to the Strominger system with
\[
\mathfrak{b}(\omega) = [|\Omega|_{\omega_{\rm CY}} \, \omega_{\rm CY}^2].
\]
The proof is by perturbation and implicit function theorem, so that the resulting solutions are of the form
\[
\omega = \omega_{\rm CY} + O(\alpha')
\]
and the Calabi-Yau metric receives higher order corrections. This result builds on earlier work of \cite{AMGF, LY05}.

The conjecture is also verified for certain balanced classes when $X$ is a $T^2$-fibration over a $K3$ surface. In this case, solutions are found in cohomology classes of the form
\[
\mathfrak{b}(\omega) = M [\omega_{K3}^2] + 2 [ \omega_{K3} \wedge i \theta \wedge \bar{\theta}], \quad M \gg 1.
\]
Here $\theta$ is a connection $(1,0)$-form on $\pi: X \rightarrow K3$. The proof given in joint work with D.H. Phong and X.-W. Zhang \cite{PPZFY} is via anomaly flow, and recovers the solutions found by Fu--Yau \cite{FY}. We will return to this example in \S \ref{dim-reduc}.

For an alternate version of the system (with auxiliary HYM tangent bundle connection $\nabla$ described in \S \ref{section-alt}), Garcia-Fernandez and Gonzalez Molina \cite{MGFM23a} define a Futaki invariant $\mathcal{F}$ which pairs with a balanced class to give zero $\langle \mathcal{F}, \mathfrak{b} \rangle = 0$ when the balanced class $\mathfrak{b}$ contains a solution to the system. It is proposed in \cite{MGFM23a} to use this invariant to find a counter-example to the conjecture, though it remains an open problem. Furthermore, it is not yet known whether there is an analogous Futaki invariant $\mathcal{F}$ for the Strominger system with Chern connection.

As for uniqueness of solutions to the Strominger system inside a given balanced class, very little is known. One of the few results in this direction is an example from \cite{FHP}, where uniqueness inside $H^4(X,\mathbb{R})$ is established assuming a warped product ansatz on a $T^4$ fibration over a Riemann surface $\Sigma_g$ of genus $g \geq 2$.

\subsection{Aeppli class} \label{section:aeppli}
Let $X$ be a complex manifold of dimension $3$ with holomorphic volume form $\Omega$ and $E \rightarrow X$ a holomorphic bundle.  Suppose $c_1^{\rm BC}(E)=0$ and $c_2^{\rm BC}(X)=c_2^{\rm BC}(E)$. Let $\omega$ be a hermitian metric on $X$ and $h$ a metric on $E$. Suppose the pair $(\omega, h)$ solves the anomaly cancellation equation
\begin{equation} \label{def-anomaly-cancellation}
i \partial \bar{\partial} \omega = \alpha' ({\rm Tr} \, R_\omega \wedge R_\omega - {\rm Tr} \, F_h \wedge F_h).
  \end{equation}
  Here $F_h = \bar{\partial} ( h^{-1} \partial h)$ and $R_\omega = \bar{\partial} (g^{-1} \partial g )$. We do not require at this point that $\omega$ is conformally balanced or $h$ is Hermitian-Yang-Mills. 

It was noticed in \cite{MGFRST} that a solution to \eqref{def-anomaly-cancellation} creates an Aeppli cohomology class. To do this, we first fix a pair of reference metrics $(\hat{\omega},\hat{h})$ and use the condition $c_2^{\rm BC}(E)=c_2^{\rm BC}(X)$ to obtain $\hat{\beta} \in \Lambda^{1,1}(X,\mathbb{R})$ satisfying
  \begin{equation} \label{c2=c2}
 {\rm Tr} \, \hat{R} \wedge \hat{R} - {\rm Tr} \, \hat{F} \wedge \hat{F} = i \partial \bar{\partial} \hat{\beta}.
  \end{equation}
The precise definition of $\hat{\beta}$ will be given below in \eqref{defn-beta}. Next, we let $R_2[g,\hat{g}] \in \Lambda^{1,1}(X,\mathbb{R})$ be the Bott-Chern-Simons secondary classes solving
  \begin{align*}
i \partial \bar{\partial} R_2[g,\hat{g}] =& {\rm Tr} \, R_\omega \wedge R_\omega - {\rm Tr} \, \hat{R} \wedge \hat{R}, \\
i \partial \bar{\partial} R_2[h,\hat{h}] =& {\rm Tr} \, F_h \wedge F_h - {\rm Tr} \, \hat{F} \wedge \hat{F}.
\end{align*}
The construction of Bott-Chern-Simons secondary classes is well-known in complex geometry (e.g. \cite{Donaldson}), but here we will use a specific construction of $R_2$ given in \eqref{defn-R2} below. Using $\hat{\beta}$ and $R_2$, we may define
  \begin{equation} \label{aeppli-class}
\mathfrak{a}(\omega,h) = \big[ \omega - \alpha' R_2[g,\hat{g}]  + \alpha' R_2[h,\hat{h}] - \alpha' \hat{\beta} \big] \in H^{1,1}_{\rm A}(X,\mathbb{R}),
  \end{equation}
  where the Aeppli cohomology is given by
  \[
H^{1,1}_{\rm A}(X,\mathbb{R}) = \frac{ (\ker \partial \bar{\partial}) \cap \Lambda^{1,1}(X,\mathbb{R}) } { {\rm Im} \, \partial \oplus \bar{\partial}}.
  \]
 We will show below that the class defined in this way does not depend on the choice of reference $(\hat{\omega},\hat{h})$. When $\alpha'=0$ and the metric $\omega$ is K\"ahler, then $\mathfrak{a}(\omega)$ is the K\"ahler class. In this sense, the Aeppli class is an $\alpha'$-corrected K\"ahler class.

Given a solution $(\omega,h)$ to \eqref{def-anomaly-cancellation}, one can look inside a fixed class for a solution to the remaining equations in the system. Concretely, this means letting
  \[
\theta \in ({\rm Im} \, \partial \oplus \bar{\partial}) \cap \Lambda^{1,1}(X,\mathbb{R})
  \]
  be an unknown variable, and substituting the ansatz
  \begin{equation} \label{aeppli-ansatz}
\tilde{\omega} = \omega + \theta  +\alpha' R_2[\tilde{g},g] - \alpha' R_2[\tilde{h},h]
  \end{equation}
 into the equations
  \[
d(|\Omega|_{\tilde{\omega}} \, \tilde{\omega}^2)=0, \quad F_{\tilde{h}} \wedge \tilde{\omega}^2=0,
  \]
  and solving for the unknowns $(\theta,\tilde{h})$. The setup is such that $(\tilde{\omega},\tilde{h})$ automatically solves the anomaly cancellation equation and $\mathfrak{a}(\tilde{\omega},\tilde{h})=\mathfrak{a}(\omega,h)$.

We note the ansatz \eqref{aeppli-ansatz} was essentially introduced by Garcia-Fernandez, Rubio, Shahbazi and Tipler \cite{MGFRST}, and there it is used to define the concept of a hermitian metric on a Bott-Chern algebroid and give a variational interpretation of the Hull-Strominger system. There is a difference betwen our setups however, and in the present case $\tilde{g}$ appears on both sides of equation \eqref{aeppli-ansatz}, unlike \cite{MGFRST} where $R_2$ is taken with respect to an auxiliary connection. Since the ansatz \eqref{aeppli-ansatz} is already itself a non-linear equation for $\tilde{g}$, further explanation is required. We give the details in the following subsection.

\subsubsection{Further details on the construction of the Aeppli class} Given two metrics $g$, $\hat{g}$, we can write
\[
{\rm Tr} \, R_g \wedge R_g  - {\rm Tr} \, R_{\hat{g}} \wedge R_{\hat{g}} = \int_0^1 {d \over dt} {\rm Tr} \, R_{g_t} \wedge R_{g_t} \, dt, \quad g_t = t g + (1-t) \hat{g},
\]
and so by \eqref{derivative-R}
\[
{\rm Tr} \, R_g \wedge R_g  - {\rm Tr} \, R_{\hat{g}} \wedge R_{\hat{g}} = i \partial \bar{\partial} \chi,
\]
with 
\[
\chi = 2 i \int_0^1  {\rm Tr} \, R_{g_t} g_t^{-1} (g-\hat{g})  \, dt.
\]
We will define $R_2[g,\hat{g}]$ using the Kodaira-Spencer \cite{KodairaSpencer} operator
\[
E = \partial \bar{\partial} (\partial \bar{\partial})^\dagger +  (\partial \bar{\partial})^\dagger \partial \bar{\partial} +  (\partial^\dagger \bar{\partial})^\dagger \partial^\dagger \bar{\partial} + \partial^\dagger \bar{\partial} (\partial^\dagger \bar{\partial})^\dagger + \bar{\partial}^\dagger \bar{\partial} + \partial^\dagger \partial,
\]
where adjoints are with respect to a reference metric $g_0$. Kodaira-Spencer \cite{KodairaSpencer} show that this is a 4th order self-adjoint elliptic operator; we refer to e.g. \cite{PopoviciBook} for an exposition of the properties of this operator. By the Fredholm alternative, there is a unique solution $\gamma \in \Lambda^{2,2} \cap (\ker E)^\perp$ solving
\[
E(\gamma) = i \partial \bar{\partial} \chi.
\]
Indeed the kernel of $E: \Lambda^{2,2}(X) \rightarrow \Lambda^{2,2}(X)$ is
\[
\ker E = \{ \eta \in \Lambda^{2,2}(X) : d \eta = 0, \ \ \bar{\partial}^\dagger \partial^\dagger \eta = 0 \} 
\]
and so $i \partial \bar{\partial} \chi$ is perpendicular to $\ker E$. Thus
\[
E(\gamma) = {\rm Tr} \, R_g \wedge R_g  - {\rm Tr} \, R_{\hat{g}} \wedge R_{\hat{g}}.
\]
Additionally, we notice that $d \gamma =0$ (see \cite{FLY}). This is because
\[
  0 = E(\gamma) - i \partial \bar{\partial} \chi
\]
can be written as
\[
0 = \partial \bar{\partial} \bigg[ \bar{\partial}^\dagger \partial^\dagger \gamma - i \chi \bigg] + \partial^\dagger \bigg[ \bar{\partial} \bar{\partial}^\dagger \partial \gamma + \partial \gamma \bigg] + \bar{\partial}^\dagger \bigg[ \partial^\dagger \partial \bar{\partial} \gamma + \partial \partial^\dagger \bar{\partial} \gamma + \bar{\partial} \gamma \bigg].
\]
The first term is orthogonal to the other two, so the sum of the last two terms is zero. Taking the inner product of the last two terms with $\gamma$ shows $\partial \gamma = 0$ and $\bar{\partial} \gamma = 0$.

Using that $d \gamma = 0$, we have
\[
E(\gamma) = \partial \bar{\partial} \bar{\partial}^\dagger \partial^\dagger \gamma = {\rm Tr} \, R_g \wedge R_g  - {\rm Tr} \, R_{\hat{g}} \wedge R_{\hat{g}},
\]
and so we define
\begin{equation} \label{defn-R2}
R_2[g,\hat{g}] = -i \bar{\partial}^\dagger \partial^\dagger E^{-1} ( {\rm Tr} \, R_g \wedge R_g  - {\rm Tr} \, R_{\hat{g}} \wedge R_{\hat{g}}).
\end{equation}
Note that we take the real part if necessary to ensure $R_2[g,\hat{g}] \in \Lambda^{1,1}(X,\mathbb{R})$. The definition of $R_2[h,\hat{h}]$ is entirely analogous. Having properly defined $R_2[h,\hat{h}]$ and $R_2[g,\hat{g}]$ in the definition of the Aeppli class \eqref{aeppli-class}, we now properly define $\hat{\beta}$ and show that the class is independent of the choice of reference metric $(\hat{\omega},\hat{h})$. For this, we start from $c_2^{\rm BC}(E)=c_2^{\rm BC}(X)$ and apply the argument detailed above to solve
\[
E(\hat{\gamma}) = \partial \bar{\partial} \bar{\partial}^\dagger \partial^\dagger \hat{\gamma} =  {\rm Tr} \, \hat{R} \wedge \hat{R} - {\rm Tr} \, \hat{F} \wedge \hat{F}
\]
and define
\begin{equation} \label{defn-beta}
\hat{\beta} = -i \bar{\partial}^\dagger \partial^\dagger \hat{\gamma}.
\end{equation}
We now verify that $\mathfrak{a}(\omega,h)$ as defined in \eqref{aeppli-class} is independent of the choice of reference $(\hat{\omega},\hat{h})$ by showing that
\[
- R_2[g,\hat{g}] + R_2[h,\hat{h}] - \hat{\beta} \in \Lambda^{1,1}(X)
\]
is independent of the choice of reference. Let $(\hat{\omega}_1,\hat{h}_1)$, $(\hat{\omega}_2,\hat{h}_2)$ be two pairs of reference metrics. Then
\[
-R_2[g,\hat{g}_1] + R_2[g,\hat{g}_2]+ R_2[h,\hat{h}_1] - R_2[h,\hat{h}_2] - \hat{\beta}_1 + \hat{\beta}_2
\]
is an expression is of the form
\[
(\partial \bar{\partial})^\dagger \Psi, \quad \Psi \in \Lambda^{2,2}(X)
\]
with $(\partial \bar{\partial}) (\partial \bar{\partial})^\dagger \Psi =0$. It follows that $ (\partial \bar{\partial})^\dagger \Psi = 0$ as desired.

We now explain why the ansatz \eqref{aeppli-ansatz} is solvable, since it is a nonlinear equation for $\tilde{g}$. First, we note that the elliptic estimate
\[
\| \gamma \|_{C^{k+4,\gamma}} \leq C \| E \gamma \|_{C^{k,\gamma}}, \quad \gamma \in (\ker E)^\perp
\]
implies that if $\hat{g}$ is a fixed smooth reference metric, then
\[
R_2[\cdot,\hat{g}]: C^{k,\gamma}(\Lambda^{1,1}) \rightarrow C^{k,\gamma}(\Lambda^{1,1}).
\]
 Using this fact, we define
  \begin{align*}
    F(\alpha',\theta,u) &= \tilde{\omega}_u - \omega - \theta - \alpha' R_2[\tilde{g}_u,g] + \alpha' R_2[\tilde{h},h] \nonumber\\
    (\tilde{g}_u)_{\alpha \bar{\beta}} &= g_{\alpha \bar{\beta}} + u_{\alpha \bar{\beta}},
\end{align*}
with
\[
F: \mathbb{R} \times C^{k,\gamma}(\Lambda^{1,1}(X,\mathbb{R})) \times  C^{k,\gamma}({\rm Herm}(X)) \rightarrow  C^{k,\gamma}(\Lambda^{1,1}(X,\mathbb{R})).
\]
Here ${\rm Herm}(X)$ denotes sections $u \in \Lambda^{1,0} \otimes \Lambda^{0,1}$ with $\overline{u_{\alpha \bar{\beta}}} = u_{\alpha \bar{\beta}}$. The equality in ansatz \eqref{aeppli-ansatz} is equivalent to $F(\alpha',\theta,u)=0$.

We note that $F(0,0,0)=0$. The derivative is
\[
D_3 F|_{(0,0,0)} = {\rm id}.
\]
By the implicit function theorem, there exists $\epsilon_1, \epsilon_2 >0$ such that for all $0<\alpha' < \epsilon_1$ and $|\theta|<\epsilon_2$, we may find a solution $u \in C^{k,\gamma}$ solving
\[
F\big(\alpha', \theta, u(\alpha',\theta) \big)= 0 .
\]
Thus for fixed $\alpha'$ small enough, the ansatz \eqref{aeppli-ansatz} is well-defined for $\theta$ taken in a suitable neighbourhood of the origin.

\subsubsection{Ellipticity}
Though the ansatz \eqref{aeppli-ansatz} allows the unknown variable $\theta$ to be in the image of $\partial \oplus \bar{\partial}$, for simplicity we take $\theta$ to be exact in this discussion. This leads to the nonlinear equation
\[
d (|\Omega|_{\omega_\theta} \, \omega^2_\theta)=0, \quad \omega_\theta = \omega + \theta + \alpha' (R_2[g_\theta,g] -R_2[\tilde{h},h] ),
\]
for an unknown $\theta \in ({\rm Im} \,d) \cap \Lambda^{1,1}(X,\mathbb{R})$. Recall that the main feature of this ansatz is that the anomaly cancellation equation is automatically implemented. We compute the linearization to check whether this can be understood as an elliptic equation. First, we linearize
\[
\delta (|\Omega|_\omega \, \omega^2) = \delta |\Omega|_\omega \, \omega^2 + 2 |\Omega|_\omega \omega \wedge \delta \omega,
\]
which using \eqref{def-norm-Omega} becomes
\begin{align*}
  \delta (|\Omega|_\omega \, \omega^2) &= |\Omega|_\omega \bigg[ -{1 \over 2} (\Lambda_\omega \delta \omega) \, \omega^2 + 2 \omega \wedge \delta \omega \bigg] \nonumber\\
  &=|\Omega|_\omega \bigg[ - \star \delta \omega + \omega \wedge \delta \omega \bigg]
\end{align*}
using the identity in complex dimension 3 (e.g. \cite{Huybrechts}) 
\[
\star \eta = - \omega \wedge \eta + {1 \over 2} (\Lambda_\omega \eta) \, \omega^2, \quad \eta \in \Lambda^{1,1}.
\]
Therefore if
\begin{align*}
  & P: C^{k,\gamma}(({\rm Im} \, d) \cap \Lambda^{1,1}) \rightarrow C^{k-1,\gamma}(\Lambda^*), \nonumber\\
  & P(\theta) = d(|\Omega|_{\omega_\theta} \omega_\theta^2),
\end{align*}
then since $d \theta = 0$ and $d \dot{\theta}=0$ we have
\begin{align*}
  \delta P(\dot{\theta}) &= d \bigg[ |\Omega|_\omega (- \star (\dot{\theta} + \alpha' \delta R_2[g_\theta,g]) + \omega \wedge (\dot{\theta} + \alpha' \delta R_2[g_\theta,g] )\bigg]  \nonumber\\
  &= -|\Omega|_\omega d \star \dot{\theta} + O(\alpha') + l.o.t.
\end{align*}
where terms without derivatives on $\dot{\theta}$ are omitted. Therefore the linearized operator is
\[
\delta P(\dot{\theta}) = |\Omega|_\omega \star (d+d^\dagger) \dot{\theta} + O(\alpha')+ l.o.t
\]
which for small $\alpha'$ is a perturbation of a 1st order elliptic operator.

From here, one could formulate a version of Yau's conjecture for Aeppli classes, in analogy to Conjecture \ref{yauconj} for balanced classes. The question remains how to formulate the stability condition for the Aeppli class $\mathfrak{a}$; in the balanced class case, the conjecture involves the notion of slope stability for $\mathfrak{b}$. In the formulation of the system with auxiliary connection $\nabla$, it is proposed in \cite{MGFRST} to formulate such a stability condition by studying the boundedness of a certain functional (the dilaton functional) which is convex along certain special paths. Uniqueness of solutions modulo automorphisms in a given Aeppli class was conjectured in \cite{MGFRST} where an analog of the geodesic equation in the space of K\"ahler metrics in a fixed K\"ahler class was introduced.

Finally, we mention Aeppli class solutions on K\"ahler backgrounds. Let $\omega_{\rm CY}$ be a Calabi-Yau metric on a K\"ahler Calabi-Yau threefold. For all $0< \alpha' < \epsilon$, we expect that there exists a solution $(\omega,h)$ of the Strominger system such that
\[
\mathfrak{a}(\omega,h) = [\omega_{\rm CY}].
\]
This result is so far only established for the system with auxiliary connection (but not yet with Chern connection): see \cite{MGFRST} for an approach via Bott-Chern algebroids and \cite{PW} for an alternate direct proof in joint work with P.-L. Wu.


\section{Equations of Motion} \label{sec:eom}
\subsection{The action functional}
The principle of least action states that to recover a physical law, we must set variations of an action to be zero. We will discuss how the Strominger system fits into this paradigm. But first, we recall two well-known examples in differential geometry.

\begin{ex}
  We start with the Einstein-Hilbert action \cite{Hilbert} on a closed manifold $X$. The action is given by
\begin{equation} \label{EinHil}
S[g_{mn}] = \int_X R \, d {\rm vol}_g,
\end{equation}
where $R$ is the scalar curvature of $g$. A classic calculation shows that imposing the variational equation $\delta S=0$ implies the Ricci-flat equation
\[
R_{mn}=0,
\]
which is the Einstein field equation from general relativity.
\end{ex}

\begin{ex}
For our second example, we consider the Riemannian metric $g$ on $X$ as fixed data, and we introduce a vector bundle $E \rightarrow X$ equipped with a metric $h$. The Yang-Mills \cite{YangMills} action is
\begin{equation} \label{YM-fun}
S[A]= \int_X |F_A|^2 \, d {\rm vol}_g,
\end{equation}
where the variable is a metric compatible connection $\nabla = d+A$, and we use the notation $F_A =dA + A \wedge A$ for the curvature of the connection $F \in \Lambda^2({\rm End} \, E)$. Imposing the variational equation $\delta S=0$, one can derive the Yang-Mills equation
\begin{equation} \label{YM-eqn}
\nabla^p  F_{p m} = 0.
\end{equation}
When $E \rightarrow X$ is a line bundle, then $F = dA \in \Lambda^2(M)$ is a 2-form representing the electromagnetic field strength and we recover H. Weyl's formulation \cite{Weyl} of the Maxwell equations in a vacuum
\begin{equation} \label{maxwell}
dF = 0, \quad \nabla^p  F_{pm} = 0.
\end{equation}
The equation $dF=0$ is the Bianchi identity. The presence of electromagnetic sources modify these equations. After introducing an electric source term, the second equation becomes $\nabla^p  F_{pm} = \rho_E$. A magnetic source term such as a magnetic monopole, which has so far never been observed in nature, would modify the first equation to $d F = \rho_M$. Though $dF \neq 0$, the equation $dF = \rho_M$ is still referred to as a Bianchi identity, though $F$ is evidently no longer the curvature of a line bundle over $X$.
\end{ex}

The examples so far have varied either the metric tensor $g_{mn}$ or the connection $\nabla = d +A$. Our next example, the action functional from heterotic supergravity, combines both Riemannian geometry and gauge theory and varies both $(g,A)$ simultaneously. There are also two more variables coming from considerations in string theory. The action for heterotic supergravity is given by \cite{BergdeRoo}
\[
  S[g,\varphi,b,A] = \int_X e^{-2 \Phi} \bigg[ R - {1 \over 12} |H|^2 + 4 |\nabla \Phi|^2 \bigg] + {\alpha' \over 8} e^{-2 \Phi} {\rm Tr} \, \bigg[ |F|^2 - |R^{\rm H}|^2 \bigg] + O(\alpha'^2)
\]
where the unknowns are:
\begin{itemize}
\item A metric tensor $g_{mn}$
\item A scalar function $\Phi: X \rightarrow \mathbb{R}$
\item A connection $\nabla = d+A$
\item A B-field $b$ with field strength $H \in \Lambda^{3}(X,\mathbb{R})$. We will discuss this object below.
\end{itemize}
The notation $R^{\rm H}$ is for the curvature tensor of the connection
\[
\Gamma^{\rm H}{}_i{}^k{}_j = \Gamma^{\rm LC}{}_i{}^k{}_j  + {1 \over 2} H_i{}^k{}_j.
\]
The B-field is a 2-form analog of the electromagnetic 1-form potential $A$ in Maxwell's equations. Its field strength is a three-form $H$, and its Bianchi identity involves an analog of a magnetic source which is produced by the structure $(g,A)$. The 3-form Bianchi identity in this setup is
\begin{equation} \label{3form-bianchi}
H = db - {\alpha' \over 4} (CS[\nabla^{\rm H}]- CS[A])
\end{equation}
where $CS[A]$ is the local Chern-Simons 3-form which satisfies $d CS[A] = {\rm Tr} \, F \wedge F$, and similarly $d CS[\nabla^{\rm H}] = {\rm Tr} \, R^{\rm H} \wedge R^{\rm H}$. The equation \eqref{3form-bianchi} has the implication that $b$ is a not a well-defined 2-form, since the local Chern-Simons 3-forms are not globally defined 3-forms. Instead, $b$ is an object which is given by local 2-form on open sets with transformation law defined such that \eqref{3form-bianchi} defines a global 3-form $H$. See \cite{CMO} for further details on B-field transformation laws. 

The action is derived as an expansion about a small parameter $\alpha'>0$. In this survey, we only consider the action $S$ up to quadratic order and denote these higher order terms by $O(\alpha'^2)$.

With this setup in place, setting the variations of the action to zero gives the equations of motion, which turn out to be (see e.g. \cite{AMP,GPT,HT})
\begin{align} \label{sugra-eom}
 & \nabla^p (e^{-2 \Phi} H_{p mn}) = O(\alpha'^2) \\
 & \hat{\nabla}^p(e^{-2 \Phi} F_{pm}) = O(\alpha') \nonumber\\ 
  &  R_{mn} + 2 \nabla_m \nabla_n \Phi - {1 \over 4} H_{mpq} H_n{}^{pq} + {\alpha' \over 4} \bigg( {\rm Tr} F_{mp} F_n{}^p - {\rm Tr} R^{\rm H}{}_{mp} R^{\rm H}{}_n{}^{p} \bigg) = O(\alpha'^2) \nonumber\\
  & R + 4 \Delta \Phi - 4 |\nabla \Phi|^2 - {1 \over 12} |H|^2 + {\alpha' \over 8} {\rm Tr} (|F|^2-|R^{\rm H}|^2)= O(\alpha'^2). \nonumber
  \end{align}
Here $\nabla$ is the Levi-Civita connection of $g$, and $\hat{\nabla}=d + \hat{\Gamma}$ is given by
\[
\hat{\Gamma}_i{}^k{}_j = \Gamma^{\rm LC}{}_i{}^k{}_j  - {1 \over 2} H_i{}^k{}_j.
\]
A theme in theoretical physics is to look for special supersymmetric solutions to the equations of motion. For non-supersymmetric solutions to these equations, see \cite{MMShah22}. Supersymmetric solutions typically involve special geometric structures, such as complex geometry or $G_2$ geometry. In this survey, we focus on special solutions arising from complex geometry.
  
\subsection{Solutions from complex geometry}

\begin{ex}
Returning to the Einstein-Hilbert functional \eqref{EinHil}, the equations of motion are $R_{mn}=0$. Let us suppose that $X$ admits a K\"ahler structure. In fact, part of K\"ahler's original motivation \cite{Kahler} for introducing K\"ahler geometry was to reduce the Ricci-flat equation to a scalar potential equation. If $X$ is a compact K\"ahler manifold, $\omega = i g_{\mu \bar{\nu}} \, dz^\mu \wedge d z^{\bar{\nu}}$ is a K\"ahler metric, and $\Omega$ is a holomorphic volume form, then the Ricci curvature of the K\"ahler metric may be written as
\begin{equation} \label{chern-conn}
R_{\mu \bar{\nu}}{}^\alpha{}_\alpha = \partial_\mu \partial_{\bar{\nu}} \log |\Omega|^2_\omega.
\end{equation}
Thus to obtain Ricci-flat metrics, we may try to solve $|\Omega|^2_{\tilde{\omega}} = e^b = const$. Looking for solutions of the form $\tilde{g}_{\mu \bar{\nu}} = g_{\mu \bar{\nu}} + \varphi_{\mu \bar{\nu}}$ where $\varphi: X \rightarrow \mathbb{R}$ is a scalar function, the Ricci-flat equation becomes the complex Monge-Amp\`ere equation
\[
\frac{\det (g_{\mu \bar{\nu}} + \varphi_{\mu \bar{\nu}}) }{\det g_{\mu \bar{\nu}}} = e^{\log |\Omega|_\omega^2 - b}.
\]
Yau's theorem \cite{Yau78} states that this equation is uniquely solvable up to shifting $\varphi$ by a constant. Thus Yau's theorem proves the existence of Ricci-flat metrics in this setting, and the resulting geometry $(X,\tilde{g},\Omega)$ is a Calabi-Yau structure.
\end{ex}

\begin{ex}
The next example is the Yang-Mills functional \eqref{YM-fun}. Let us suppose that $E \rightarrow X$ is equiped with additional structure: suppose $X$ is a complex manifold with fixed K\"ahler metric $\omega$. Then the equations
\begin{equation} \label{HYM-eqn}
F^{0,2}=0, \quad g^{\mu \bar{\nu}} F_{\mu \bar{\nu}} = 0
\end{equation}
are well-known to be special solutions to the Yang-Mills equations \cite{AHS}. Indeed, starting from the Bianchi identity $\bar{\partial}_A F =0$ and converting derivatives to the Chern connection, we obtain
\[
\nabla_{\bar{\beta}} F_{\mu \bar{\nu}} - \nabla_{\bar{\nu}} F_{\mu \bar{\beta}} =0,
\]
and upon taking the trace with $g^{\mu \bar{\nu}}$ we obtain $\nabla^\mu F_{\mu \bar{\nu}} = 0$ which is the Yang-Mills equation \eqref{YM-eqn}. Equation \eqref{HYM-eqn} is the Hermitian-Yang-Mills equation, and the existence of solutions is equivalent to a notion of stability in algebraic geometry by the Donaldson-Uhlenbeck-Yau theorem \cite{Donaldson, UY}.
\end{ex}

We return to our main topic of heterotic supergravity. For this, let $X$ be a complex manifold of dimension 3 with hermitian metric $\omega$ and holomorphic volume form $\Omega$. Let
\[
H = i (\partial - \bar{\partial})\omega, \quad \Phi = - {1 \over 2} \log |\Omega|_\omega.
\]
Note that if $\omega$ is a K\"ahler Ricci-flat metric, then $H= 0$ and $\Phi$ is constant. The equations of motion \eqref{sugra-eom} can be solved if $E=T^{1,0}X$ by taking $A$ to be the Chern connection of the K\"ahler Ricci-flat metric. These are the solutions found by Candelas-Horowitz-Strominger-Witten \cite{CHSW} which introduced Calabi-Yau geometry into string theory.

More generally, given a vector bundle $E \rightarrow X$, we may look for solutions to the system without setting $H=0$, so that $\omega$ is non-K\"ahler. The supersymmetric equations \cite{Strominger,Hull} for the heterotic system are
\begin{align} \label{susy1}
d(|\Omega|_\omega \, \omega^2)= 0, \quad F \wedge \omega^2 = 0, \quad F=\bar{\partial}( h^{-1}\partial h) \\
\label{susy2}
i \partial \bar{\partial} \omega = {\alpha' \over 8} (  {\rm Tr} \, R^{\rm H} \wedge R^{\rm H} - {\rm Tr} \, F \wedge F) + O(\alpha'^2).
\end{align}

It is well-known in the string theory literature that these equations are special solutions to the equations of motion; see e.g. \cite{AMP,CCDL} or \cite{Ivanov} with different connections on the tangent bundle. For completeness, we give a direct proof of this fact, namely that solutions to the equations \eqref{susy1}, \eqref{susy2} satisfy the equations of motion \eqref{sugra-eom}.

\subsection{Divergence equations}
We start by showing that conformally balanced metrics solve the following Yang-Mills type equation for the 3-form field strength $H$:
\begin{equation} \label{div-Ho}
(\nabla^{\rm LC}){}^p (e^{-2\Phi} H_{mnp}) = 0.
\end{equation}
We recall our convention that $m,n,p$ denote real indices while $\mu,\nu,\alpha$ denote complex holomorphic indices. Before computing the divergence of $H$, we rewrite the conformally balanced condition in a different form which will be frequently used in the upcoming calculations. A direct calculation (e.g. \cite{LY05,PPZ18CAG,cetraro}) shows that $d(|\Omega|_\omega \, \omega^2)=0$ is equivalent to
\[
H_\mu{}^\mu{}_\alpha = \partial_\alpha \log |\Omega|_\omega, \quad H_\mu{}^\mu{}_{\bar{\beta}} = - \partial_{\bar{\beta}} \log |\Omega|_\omega.
\]
Writing this constraint in terms of the scalar function $\Phi=-{1 \over 2} \log |\Omega|_\omega$ gives
\begin{equation} \label{confbal-dilaton}
H_\mu{}^\mu{}_\alpha = -2 \Phi_\alpha, \quad H_\mu{}^\mu{}_{\bar{\beta}} = 2 \Phi_{\bar{\beta}}.
\end{equation}
We now show that \eqref{confbal-dilaton} implies \eqref{div-Ho}. We start with
\begin{equation} \label{div-H0}
(\nabla^{\rm LC}){}^m  H_{m \alpha \beta} = (\nabla^{\rm LC}){}^\lambda  H_{\lambda \alpha \beta} +  (\nabla^{\rm LC}){}^{\bar{\lambda}}  H_{\bar{\lambda} \alpha \beta}.
\end{equation}
By the expression for the Levi-Civita connection in holomorphic coordinates, we have
\begin{align}
 (\nabla^{\rm LC}){}^\lambda H_{\lambda \alpha \beta} &= g^{\lambda \bar{\nu}} (\partial_{\bar{\nu}} H_{\lambda \alpha \beta} - \Gamma_{\bar{\nu}}{}^n{}_\lambda H_{n \alpha \beta} - \Gamma_{\bar{\nu}}{}^n{}_\alpha H_{\lambda n \beta} - \Gamma_{\bar{\nu}}{}^n{}_\beta H_{\lambda \alpha n}) \nonumber\\
 \label{div-H00} &= - {1 \over 2} H^{\lambda \bar{\sigma}}{}_\lambda H_{\bar{\sigma} \alpha \beta}, 
\end{align}
using  \eqref{LC2C} and $H_{\lambda \alpha \beta} = 0$. Next, we compute the second term in \eqref{div-H0} by converting to the Chern connection via \eqref{LC2C}.
\begin{equation} \label{div-H1}
(\nabla^{\rm LC}){}^{\bar{\lambda}}  H_{\bar{\lambda} \alpha \beta} = \nabla^{\rm Ch}{}^{\bar{\lambda}} H_{\bar{\lambda} \alpha \beta} - {1 \over 2} H_\mu{}^{\mu \bar{\lambda}} H_{\bar{\lambda} \alpha \beta}.
\end{equation}
 Next, we note the general identity
\begin{equation} \label{H-mysteryswitch}
\nabla^{\rm Ch}{}^{\bar{\lambda}} H_{\bar{\lambda} \alpha \beta} = - \partial_\alpha H_\mu{}^\mu{}_\beta + \partial_\beta H_\mu{}^\mu{}_\alpha,
\end{equation}
which follows from cancelling mixed partial derivatives in the expression $H_{\mu \bar{\nu} \beta} = \partial_\mu g_{\bar{\nu} \beta} - \partial_\beta g_{\bar{\nu} \mu}$. We now substitute the conformally balanced condition \eqref{confbal-dilaton} which implies that conformally balanced metrics satisfy
\begin{equation} \label{div-H-Ch}
\nabla^{\rm Ch}{}^{\bar{\lambda}} H_{\bar{\lambda} \alpha \beta} = 0.
\end{equation}
Adding \eqref{div-H00} and \eqref{div-H1} and substituting \eqref{div-H-Ch}, \eqref{confbal-dilaton} gives the identity
\[
(\nabla^{\rm LC}){}^m  H_{m \alpha \beta} = 2 \Phi^{\bar{\lambda}} H_{\bar{\lambda} \alpha \beta},
\]
which holds for conformally balanced metrics. Therefore
\[
(\nabla^{\rm LC}){}^m (e^{-2\Phi} H_{m \alpha \beta}) = 0.
\]
It remains to show \eqref{div-Ho} when $m,n$ are mixed barred/unbarred holomorphic indices. Convertion of the Levi-Civita connection to the Chern connection together with the conformally balanced identity \eqref{confbal-dilaton} gives the identity
\begin{equation} \label{div-H2}
\nabla^{\rm LC}{}^{\bar{\lambda}} H_{\bar{\lambda} \alpha \bar{\beta}}  +  \nabla^{\rm LC}{}^{\lambda} H_{\lambda \alpha \bar{\beta}} = \nabla^{\rm Ch}{}^{\bar{\lambda}} H_{\bar{\lambda} \alpha \bar{\beta}}  +  \nabla^{\rm Ch}{}^{\lambda} H_{\lambda \alpha \bar{\beta}}  + 2 \Phi^m H_{m \alpha \bar{\beta}}.
\end{equation}
In a similar way to \eqref{H-mysteryswitch}, there is the following general identity in Hermitian geometry
\[
 \nabla^{\rm Ch}{}^{\bar{\lambda}} H_{\bar{\lambda} \alpha \bar{\beta}}  +  \nabla^{\rm Ch}{}^{\lambda} H_{\lambda \alpha \bar{\beta}} = - \partial_\alpha H_\mu{}^\mu{}_{\bar{\beta}} - \partial_{\bar{\beta}} H_\mu{}^\mu{}_\alpha.
\]
This vanishes by the conformally balanced identity \eqref{confbal-dilaton}. Substituting this into \eqref{div-H2} gives
\[
(\nabla^{\rm LC}){}^m (e^{-2\Phi} H_{m \alpha \bar{\beta} }) = 0.
\]
Thus \eqref{div-Ho} is verified for conformally balanced metrics.
\smallskip
\par One of the four equations of motions is then satisfied. We move on to verify the non-K\"ahler Yang-Mills equation
\begin{equation} \label{div-F}
\hat{\nabla}^m (e^{-2\Phi} F_{mn}) = 0.
\end{equation}
This equation holds for Hermitian-Yang-Mills connection over a conformally balanced manifold. We note in particular that HYM connections over non-K\"ahler geometries are do not solve the Yang-Mills equation \eqref{YM-eqn}, as the classic Yang-Mills equation is not the equation required by heterotic supergravity when $H \neq 0$ and $\Phi \neq const$. Converting the Strominger-Bismut connection $\hat{\nabla}$ to the Chern connection via \eqref{Str-Bis} gives
\[
\hat{\nabla}_\kappa F^{\kappa \bar{\mu}} = \nabla^{\rm Ch}{}_\kappa F^{\kappa \bar{\mu}} + H^\kappa{}_{\kappa \lambda} F^{\lambda \bar{\mu}} + H^{\bar{\mu}}{}_{\kappa \bar{\lambda}} F^{\kappa \bar{\lambda}}.
\]
The Bianchi identity $d_AF=0$ using the Chern connection can be written as
\[
\nabla^{\rm Ch}{}_{\bar{\kappa}} F_{\mu \bar{\nu}} - \nabla^{\rm Ch}{}_{\bar{\nu}} F_{\mu \bar{\kappa}} = H_{\bar{\kappa} \bar{\nu}}{}^{\bar{\lambda}} F_{\mu \bar{\lambda}}.
\]
Using $g^{\mu \bar{\nu}} F_{\mu \bar{\nu}} = 0$, this implies
\[
\hat{\nabla}_\kappa F^{\kappa \bar{\mu}} = H^\kappa{}_{\kappa \lambda} F^{\lambda \bar{\mu}}.
\]
Using  the conformally balanced condition \eqref{confbal-dilaton} proves \eqref{div-F}. 

\subsection{Einstein equation}
We now verify the third equation in \eqref{sugra-eom}, which is the Einstein equation in heterotic string theory. The first step is to compare the Riemannian Ricci curvature with the curvature of the Chern connection. A long direct calculation converting \eqref{curvature}, \eqref{ricci} into the Chern connection via \eqref{LC2C} leads to (see e.g. \cite{AMP} for details)
\begin{align} \label{Ricci}
  R_{\alpha \beta} &= {1 \over 2} (\nabla^{\rm Ch}_\alpha H_\mu{}^\mu{}_\beta + \nabla^{\rm Ch}_\beta H_\mu{}^\mu{}_\alpha) + \frac{1}{2} H_{\alpha \mu \bar{\nu}} H^{\mu \bar{\nu}}{}_\beta \\
  R_{\alpha \bar{\beta}} &= R^{\rm Ch}{}_\mu{}^\mu{}_{\bar{\beta} \alpha} - {1 \over 2}(\nabla^{\rm Ch}_\alpha H_\mu{}^\mu{}_{\bar{\beta}} - \nabla^{\rm Ch}_{\bar{\beta}} H_\mu{}^\mu{}_\alpha) \nonumber\\
 & \ \ -{1 \over 2} H_{\alpha \mu \bar{\nu}} H^{\mu \bar{\nu}}{}_{\bar{\beta}} + {1 \over 4} H_{\alpha \bar{\mu} \bar{\nu}} H^{\bar{\mu} \bar{\nu}}{}_{\bar{\beta}} - {1 \over 2} H_\alpha{}^m{}_\beta H_\mu{}^\mu{}_m. \nonumber
\end{align}
If $H=0$, we recover the formula from K\"ahler geometry which identifies the Ricci curvature $R_{\alpha \bar{\beta}}$ with the trace of the Chern curvature. Next, we use the following identity for the Chern connection (see e.g. Appendix B in \cite{MPS})
\begin{align} \label{R-dH}
  R^{\rm Ch}{}_\mu{}^\mu{}_{\bar{\beta} \alpha} &= R^{\rm Ch}{}_{\alpha \bar{\beta}}{}^\mu{}_\mu - (i \partial \bar{\partial} \omega)_\mu{}^\mu{}_{\bar{\beta} \alpha} \\
  & \ \ + \nabla^{\rm Ch}_\alpha H_\mu{}^\mu{}_{\bar{\beta}} - \nabla^{\rm Ch}_{\bar{\beta}} H_\mu{}^\mu{}_\alpha + H_{\alpha \mu \bar{\nu}} H^{\mu \bar{\nu}}{}_{\bar{\beta}}. \nonumber
\end{align}
When $\omega$ is K\"ahler, this is simply the symmetry $R_\mu{}^\mu{}_{\bar{\beta} \alpha} = R_{\alpha \bar{\beta}}{}^\mu{}_\mu$ which is obstructed by non-K\"ahler terms in general.

The identities \eqref{Ricci} and \eqref{R-dH} hold for general hermitian geometries. We now specialize to conformally balanced metrics.
Substituting the conformally balanced constraint \eqref{confbal-dilaton} together with \eqref{chern-conn} into \eqref{R-dH} gives the identity
\begin{equation} \label{Chern-Ricci-confbal}
R^{\rm Ch}{}_\mu{}^\mu{}_{\bar{\beta} \alpha} = - (i \partial \bar{\partial} \omega)_\mu{}^\mu{}_{\bar{\beta} \alpha} + H_{\alpha \mu \bar{\nu}} H^{\mu \bar{\nu}}{}_{\bar{\beta}}
\end{equation}
which holds for conformally balanced metrics. Substituting \eqref{Chern-Ricci-confbal} into \eqref{Ricci} gives \cite{AMP}
\begin{align} \label{eom-ddbar}
  & R_{\alpha \bar{\beta}} + 2 \nabla^{\rm LC}{}_\alpha \nabla^{\rm LC}{}_{\bar{\beta}} \Phi - {1 \over 4} H_{\alpha mn} H_{\bar{\beta}}{}^{mn} = (- i \partial \bar{\partial} \omega)_\mu{}^\mu{}_{\bar{\beta} \alpha} \\
  & R_{\alpha \beta} + 2 \nabla^{\rm LC}{}_\alpha \nabla^{\rm LC}{}_{\beta} \Phi - {1 \over 4} H_{\alpha mn} H_\beta{}^{mn} = 0. \nonumber
\end{align}
Recall our convention is such that Greek indices $\alpha, \beta$ are holomorphic indices, while Roman indices $m,n$ are real. We now move on to the $\alpha'$-terms in \eqref{sugra-eom}. Our conventions are such that
\[
({\rm Tr} \, F \wedge F)_{\mu \bar{\nu} \alpha \bar{\beta}} = 2 {\rm Tr} \, F_{\mu \bar{\nu}} F_{\alpha \bar{\beta}} - 2 {\rm Tr} \, F_{\mu \bar{\beta}} F_{\alpha \bar{\nu}}.
\]
Therefore assuming
\[
F^{0,2}=0, \quad g^{\mu \bar{\nu}} F_{\mu \bar{\nu}}= 0 , \quad g^{\mu \bar{\nu}} R^{\rm H}{}_{\mu \bar{\nu}} = O(\alpha'), \quad (R^{\rm H}){}^{0,2} = O(\alpha')
\]
we obtain
\begin{align} \label{Einstein-alpha}
& - \alpha' \bigg( {\rm Tr} \, F_{\alpha p} F_{\bar{\beta}}{}^p - {\rm Tr} \, R^{\rm H}{}_{\alpha p} R^{\rm H}{}_{\bar{\beta}}{}^{p} \bigg) \nonumber\\
  &= {\alpha' \over 2} ({\rm Tr} \, F \wedge F - {\rm Tr} \, R^{\rm H} \wedge R^{\rm H})_\mu{}^\mu{}_{\bar{\beta} \alpha} + O(\alpha'^2).
\end{align}
Here the trace in ${\rm Tr} \, R \wedge R$ is over all real indices (over $T_{\mathbb{C}} X$). Combining \eqref{eom-ddbar} and \eqref{Einstein-alpha} and applying the anomaly relation \eqref{susy2} verifies \eqref{sugra-eom}. 

We now explain $g^{\mu \bar{\nu}} R ^{\rm H}{}_{\mu \bar{\nu}} = O(\alpha')$. Direct calculation of $R ^{\rm H}{}_{\mu \bar{\nu}}$, as noted in the appendix \eqref{RH11}, gives
\[
g^{\mu \bar{\nu}} R^{\rm H}{}_{\mu \bar{\nu} \bar{\beta} \alpha}= R^{\rm Ch}{}^\mu{}_{\mu \bar{\beta} \alpha} - H_\alpha{}^{\mu \bar{\sigma}} H_{\mu \bar{\sigma} \bar{\beta}} .
\]
Conformally balanced metrics satisfy the identity \eqref{Chern-Ricci-confbal}, hence
\[
g^{\mu \bar{\nu}} R^{\rm H}{}_{\mu \bar{\nu} \bar{\beta} \alpha} = - (i \partial \bar{\partial} \omega)_\mu{}^\mu{}_{\bar{\beta} \alpha} = O(\alpha')
\]
by the anomaly relation \eqref{susy2}. We also note that
\[
g^{\mu \bar{\nu}} R^{\rm H}{}_{\mu \bar{\nu} \alpha \beta} = \nabla^{\rm Ch}{}^{\bar{\nu}} H_{\bar{\nu} \alpha \beta} = 0
\]
by \eqref{div-H-Ch}. We also need to verify $(R^{\rm H}){}^{0,2} = O(\alpha')$. As noted in appendix \eqref{RH20}, the only non-zero component of $R_{\mu \nu m n}$ is present when both $m$ and $n$ are barred holomorphic indices, and the expression is
\[
R^{\rm H}{}_{\mu \nu \bar{\alpha} \bar{\beta}} = - (i \partial \bar{\partial} \omega)_{\mu \nu \bar{\alpha} \bar{\beta}} = O(\alpha')
\]
as required.

\subsection{Dilaton equation}
It remains to verify the last equation in \eqref{sugra-eom}. The terms can be regrouped as
\begin{align*}
 &  \bigg( R + 2 \Delta \Phi - {1 \over 4} |H|^2 + {\alpha' \over 4}({\rm Tr} \, F_{mn} F^{mn} - {\rm Tr} \, R^{\rm H}{}_{mn} R^{\rm H}{}^{mn}) \bigg) \nonumber\\
  &+ \bigg(2 \Delta \Phi + {1 \over 6} |H|^2 - 4 |\nabla \Phi|^2 - {\alpha' \over 8} ({\rm Tr} \, F_{mn} F^{mn} - {\rm Tr} \, R^{\rm H}{}_{mn} R^{\rm H}{}^{mn}) \bigg).
\end{align*}
The first grouping is zero to $O(\alpha'^2)$ by taking the contraction of the heterotic Einstein equation, which is
\begin{equation} \label{Trace-Ein}
R + 2 \Delta \Phi - {1 \over 4} |H|^2 + {\alpha' \over 4} {\rm Tr} \, (|F|^2-|R^{\rm H}|^2) = O(\alpha'^2).
\end{equation}
We are left with the second grouping, for which we use identities from complex geometry. Since $\Phi = - {1 \over 2} \log |\Omega|_\omega$, we have
\[
g^{\alpha \bar{\beta}} \Phi_{\alpha \bar{\beta}} = - {1 \over 4} R^{\rm Ch}{}_{\alpha}{}^\alpha{}^\mu{}_\mu.
\]
We next use \eqref{Chern-Ricci-confbal}, which holds for conformally balanced metrics and implies
\begin{equation} \label{dilaton1}
g^{\alpha \bar{\beta}} \Phi_{\alpha \bar{\beta}} = {1 \over 4} (i \partial \bar{\partial} \omega)_{\alpha}{}^\alpha{}^\mu{}_\mu - {1 \over 4} H_{\alpha \mu \bar{\nu}} H^{\mu \bar{\nu} \alpha}.
\end{equation}
Direct conversion of covariant derivatives with respect to the Levi-Civita connection gives
\begin{equation} \label{dilaton2}
g^{\alpha \bar{\beta}} \nabla^{\rm LC}{}_\alpha \nabla^{\rm LC}{}_{\bar{\beta}} \Phi = g^{\alpha \bar{\beta}} \Phi_{\alpha \bar{\beta}} + {g^{\alpha \bar{\beta}} \over 2} (-H_\alpha{}^\mu{}_{\bar{\beta}} \Phi_\mu + H_\alpha{}^{\bar{\mu}}{}_{\bar{\beta}} \Phi_{\bar{\mu}}).
\end{equation}
We can now combine \eqref{dilaton1} and \eqref{dilaton2}, apply \eqref{confbal-dilaton}, and convert all sums to real coordinates (so that e.g. $\Delta \Phi = 2 g^{\alpha \bar{\beta}} \nabla^{\rm LC}{}_\alpha \nabla^{\rm LC}{}_{\bar{\beta}} \Phi$ and $|H|^2 = 6 H_{\alpha \mu \bar{\nu}} H^{\mu \bar{\nu} \alpha}$) to deduce
\begin{equation} \label{dilaton3}
\Delta \Phi = - {1 \over 12} |H|^2 + 2 |\nabla \Phi|^2 + {1 \over 2} (i \partial \bar{\partial} \omega)_\alpha{}^{\alpha \mu}{}_\mu.
\end{equation} 
Lastly, we note that \eqref{Einstein-alpha} implies
\[
\alpha' ({\rm Tr} \, F_{mn} F^{mn} - {\rm Tr} \, R^{\rm H}{}_{mn} R^{\rm H}{}^{mn}) = -\alpha' ({\rm Tr} \, F \wedge F - {\rm Tr} \, R \wedge R)_\alpha{}^{\alpha \mu}{}_\mu + O(\alpha'^2).
\]
Combining this identity, \eqref{dilaton3}, and the heterotic Bianchi identity \eqref{susy2} implies that the second grouping is zero to $O(\alpha'^2)$.

\subsection{Strominger system and equations of motion} \label{Strom-eom}
In the previous sections we assumed the equations \eqref{susy1}, \eqref{susy2} and presented the well-known result from the string theory literature that these complex geometries satisfy the equations of motion \eqref{sugra-eom}. The difference with the setup considered in the current survey is that we use the Chern connection to compute ${\rm Tr} \, R \wedge R$, and we demand the exact equality
\begin{equation} \label{bianchi-chern}
i \partial \bar{\partial} \omega = {\alpha' \over 8} ( {\rm Tr} \, R^{\rm Ch} \wedge R^{\rm Ch} - {\rm Tr} \, F \wedge F).
\end{equation}
The question then becomes: when does a solution to the PDE \eqref{bianchi-chern} give a solution to the physical $\alpha'$ expansion \eqref{susy2}? Since \eqref{susy2} solves the equations of motion, this would give a condition under which solutions to the Strominger system with Chern connection solve the equations of motion. For this, we use the identities from the appendix \eqref{RH11} to convert the Chern curvature to the Hull curvature and obtain
\begin{align*}
  & ({\rm Tr} \, R^{\rm Ch} \wedge R^{\rm Ch}) = ({\rm Tr} \, R^{\rm H} \wedge R^{\rm H})^{2,2} + R^{\rm Ch} * H^2 + H^4 + (\nabla H)^2 \\
 & ({\rm Tr} \, R^{\rm H} \wedge R^{\rm H})^{3,1} = \nabla H * (dH) \\
  & ({\rm Tr} \, R^{\rm H} \wedge R^{\rm H})^{4,0} = 0.
\end{align*}

To solve the physical $\alpha'$ expansion and hence the equations of motion, we require a solution to the Strominger system with Chern connection together with an estimate.

\begin{prop} \label{chern-eom}
  Let $X$ be a compact complex manifold of dimension 3 with holomorphic volume form $\Omega$. Let $E \rightarrow X$ be a holomorphic vector bundle. Let $(\omega,h) = (\omega_{\alpha'},h_{\alpha'})$ be a sequence of solutions to
  \begin{align}
&    d(|\Omega|_\omega \, \omega^2) = 0, \quad F_h \wedge \omega^2 = 0 \nonumber\\
& i \partial \bar{\partial} \omega = {\alpha' \over 8} ( {\rm Tr} \, R^{\rm Ch} \wedge R^{\rm Ch} - {\rm Tr} \, F \wedge F). \nonumber
\end{align}
along a sequence $\alpha' \rightarrow 0$, together with estimates
\begin{equation} \label{H-est}
|H| + |\nabla H| \leq C \sqrt{\alpha'}, \quad |F| + |R^{\rm Ch}| \leq C.
\end{equation}
Then the 4-form condition \eqref{susy2} with Hull connection is satisfied to linear order in $\alpha'$, and consequently the solutions $(\omega,h)$ satisfy the equations of motion to linear order in $\alpha'$ in the action $S$.
\end{prop}

In other words, solutions to the Strominger system with $H = O(\sqrt{\alpha'})$ are solutions to the equations of motion of the functional $S$. This fact can also be found in the string theory literature e.g. \cite{McOSva,MMS}. This indicates that physical solutions to the Strominger system should satisfy extra PDE estimates. This is reasonable as nonlinear PDE are typically solved in an appropriate regime; for example when verifying the ellipticity of the Strominger system, an estimate of the form $|\alpha' R^{\rm Ch}|< 1$ is sufficient, and this is implied by estimate \eqref{H-est}.

We give two examples of solutions satisfying the estimate \eqref{H-est}:

\begin{ex}
Perturbations about a K\"ahler solution. The solutions obtained by Li--Yau \cite{LY05} (see Andreas--Garcia-Fernandez \cite{AMGF} for a generalization of this result, and see \cite{CPY-IFT} for an alternate proof) are smooth perturbations by the implicit function theorem about an $\alpha'=0$ K\"ahler Calabi-Yau solution. This means that $g \rightarrow g_{\rm CY}$ smoothly as $\alpha' \rightarrow 0$, and so $|H| + |\nabla H| \leq C \alpha'$ by the implicit function theorem. For a discussion of the equations of motion near a K\"ahler solution in the string theory literature, we refer to \cite{MMS,Witten2}.
\end{ex}

\begin{ex}
Solutions on the total space of a holomorphic $T^2$ fibration over a $K3$ surface: $\pi: X \rightarrow K3$. The fact that these solve the heterotic equations of motion is explained in \cite{MMS} and \cite{AMP}. We present here an outline of the main ideas from the perspective of the Fu--Yau \cite{FY} solutions to the system with Chern connection.

The metric ansatz in this geometric setup is parametrized by a scalar function $u$ on the base $K3$:
\[
\omega = e^u \omega_{\rm K3}+ a i \theta \wedge \bar{\theta},
\]
where $\theta = \theta_1+i \theta_2$ is a connection $(1,0)$-form with $d \theta = \omega_1+i \omega_2$ and $\omega_i \in 2 \pi H^2(K3,\mathbb{Z})$ is a pair of anti-self-dual forms with respect to the K\"ahler Ricci-flat metric $\omega_{\rm K3}$. The vector bundle metric is of the form $(\pi^* E_{K3},\pi^* h_{K3})$ where $E_{K3} \rightarrow K3$ is a stable bundle over the $K3$ with Hermitian-Yang-Mills metric $h_{K3}$ solving
\[
F_{h_{K3}} \wedge \omega_{\rm K3} = 0.
\]
The ansatz can be substituted into the Strominger system with Chern connection, and it leads to a nonlinear PDE for $u$ which is solved by Fu and Yau \cite{FY}.

As explained in the alternate approaches of \cite{MMS,MGF-T}, by setting $a = \alpha'$ there are examples of vector bundles $E_{K3}$ for which the system is solvable for each $\alpha'>0$. In the language of \cite{FY,PPZFY}, the integrability condition $\int \mu = 0$ is satisfied along $\alpha' \rightarrow 0$.

At each value of $\alpha'$, we may start the anomaly flow with $t=0$ initial data $u(0) = \log M$ where $M \gg 1$ is a constant. The PDE estimates on the scalar function $u$ derived in \cite{PPZFY} give convergence to a solution to the Strominger system as $t \rightarrow \infty$ with bounds of the form
\begin{equation} \label{FY-est}
{1 \over CM} \leq e^{-u} \leq {C \over M}, \quad |D^k u|_{g_{\rm K3}} \leq C_k, \quad k \geq 1,
\end{equation}
where $C$ is uniform over $0<\alpha' \leq 1$ and $0 < a \leq 1$. Applying this estimate to
\[
H^{1,2} = i \bar{\partial} \omega = i e^u \bar{\partial} u \wedge \omega_{\rm K3} - a (\omega_1+i \omega_2) \wedge \bar{\theta},
\]
gives
\[
|H| + |\nabla H| \leq C M^{-1/2}.
\]
This suggests to choose the scale $M= \alpha'^{-1}$. The Chern curvature $R$ is computed in \cite{FY} (see also \cite{cetraro} for another exposition), and the estimates \eqref{FY-est} imply that $|R^{\rm Ch}| \leq C$. The connection on the gauge bundle is fixed in this ansatz, so $|F| \leq CM^{-1} \leq C$. Thus we may solve the Strominger system along a sequence $\alpha' \rightarrow 0$ for which the estimates \eqref{H-est} are satisfied. 
\end{ex}

\section{Anomaly Flow}
\subsection{Definitions}
Let $(M,g)$ be a compact Riemannian manifold. Hamilton's Ricci flow equation \cite{Hamilton, HamiltonOmnibus} is
\[
\partial_t g_{ij} = - 2 R_{ij}.
\]
This nonlinear PDE for the metric tensor is a differential geometric analog of the heat equation. This analogy can be understood in harmonic coordinates $\{ x^i \}$, where when viewing $g_{ij}(x)$ as a local function in these coordinates we have the identity
\[
-2 R_{ij} = \Delta_g g_{ij}(x) + \mathcal{O}(g,Dg).
\]
As a sort of heat flow for the metric tensor, the Ricci flow is intuitively expected to improve the geometry of the underlying manifold.

If the manifold $(M,g)$ is equipped with additional structure, there are versions of the Ricci flow equation adapted to various special geometries. 

\begin{ex}
If the manifold $M$ admits a complex structure and K\"ahler metric $g$, then the Ricci flow starting from a K\"ahler metric preserves the K\"ahler condition and is referred to as the K\"ahler-Ricci flow:
\[
\partial_t \omega = - {\rm Ric}(\omega), \quad d \omega_0 = 0.
\]
The study of K\"ahler-Ricci flow was initiated by Cao \cite{Cao}. We refer to \cite{KRF, KAWA} for modern developments and applications in K\"ahler geometry. 
\end{ex}

\begin{ex}
There is an analog of Ricci flow for holomorphic vector bundles introduced by Donaldson \cite{Donaldson}. Let $h$ be a hermitian metric on a holomorphic vector bundle $E \rightarrow M$ with $c_1(E)=0$ over a complex manifold with K\"ahler metric $\omega$. Let $F_h = \bar{\partial}(h^{-1} \partial h)$ denote the curvature tensor of the metric $h$. The Donaldson heat flow is:
\begin{equation} \label{donaldson-heat}
h^{-1} \partial_t h = - i \Lambda_\omega F_h.
\end{equation}
Here $i \Lambda_\omega F_h \in \Gamma({\rm End} \, E)$ with $i \Lambda_\omega F_h = g^{\mu \bar{\nu}} F_{\mu \bar{\nu}}$. We note that this flow is not the K\"ahler-Ricci flow when $E=T^{1,0}X$ since the background metric appearing in $\Lambda_\omega$ is fixed, but it is an example of a geometric flow involving the trace of curvature.  
\end{ex}

Having discussed various aspects of the Strominger system in the first part of this survey, we now examine a Ricci flow-type setup in this special geometry. There is a parabolic version of the Strominger system which will flow the metric by the Ricci curvature plus extra corrections. The anomaly flow, introduced in joint work with D.H. Phong and X.-W. Zhang \cite{PPZ18}, is the flow of a pair $(\omega(t),h(t))$ determined by the equation
\begin{align} \label{AF}
 \partial_t (|\Omega|_{\omega(t)} \, \omega(t)^2) &= i \partial \bar{\partial} \omega(t) - \alpha' (  {\rm Tr} \, R_{\omega(t)} \wedge R_{\omega(t)} - {\rm Tr} \, F_{h(t)} \wedge F_{h(t)}), \nonumber\\
d (|\Omega|_{\omega_0} \, \omega_0^2) &= 0.
\end{align}
We recall our setup: $X$ is a compact complex manifold of dimension 3 with fixed holomorphic volume form $\Omega$, and $E \rightarrow X$ is a holomorphic vector bundle with $c_1(E)=0$ and $c_2^{\rm BC}(E)=c_2^{\rm BC}(X)$. The metric $h(t)$ on the bundle flows by the Donaldson heat flow \eqref{donaldson-heat} but with varying background hermitian metric $\omega(t)$.

Our first remark is on connecting the anomaly flow with the discussion on cohomology classes from \S \ref{section:coho}. We presented there two cohomology classes associated to the Strominger system: the Aeppli class $\mathfrak{a} \in H^{1,1}_{\rm A}(X,\mathbb{R})$ and the balanced class $\mathfrak{b} \in H^{2,2}_{\rm BC}(X,\mathbb{R})$. The anomaly flow takes the balanced class point of view. Namely, the balanced condition is preserved for all time
\[
d (|\Omega|_{\omega(t)} \, \omega(t)^2) = 0,
\]
and the solution remains in the balanced class of the initial metric:
\[
|\Omega|_{\omega(t)} \, \omega(t)^2 \in [ |\Omega|_{\omega_0} \, \omega_0^2 ] = \mathfrak{b}.
\]
In other words, the anomaly flow takes a conformally balanced metric as initial data, and seeks a solution to the Strominger system in the fixed balanced class given by the initial metric. In the context of \S \ref{section:coho}, we wonder if there is a natural flow seeking a Strominger system representative in a fixed Aeppli class. 

The anomaly flow has a short-time existence and uniqueness property. It is shown in \cite{PPZ18} that \eqref{AF}, though it is given as a flow of $(2,2)$ forms, does uniquely specify a flow for the metric $\omega(t)$, and short-time existence is established by linearizing the equation and proving the invertibility of the symbol. However, the flow \eqref{AF} is a fully nonlinear PDE for the metric due to the term ${\rm Tr} \, R_{\omega(t)} \wedge R_{\omega(t)}$, and so short-time existence should not be expected for arbitrary initial data; fully nonlinear equations generally have an open region of admissible solutions where the linearization is elliptic. In this case, the flow is parabolic once the terms of order $\alpha'$ are small.
\begin{thm}
\cite{PPZ18} If $\omega_0$ is an initial conformally balanced metric with $|\alpha' R_{g_0}|< {1 \over 2}$, then there exists $T>0$ such that the anomaly flow \eqref{AF} admits a unique solution on $[0,T)$ with $\omega(0) = \omega_0$.
\end{thm}
For fixed $\alpha'>0$, the condition $|\alpha' R_{g_0}|< {1 \over 2}$ can be achieved by an initial metric after scaling $g \mapsto \lambda g$. There are examples (see e.g. \cite{FHP2, PPZFY}) of the condition $|\alpha' R| \ll 1$ being preserved along the flow for all $t \in [0,\infty)$, however the general theory and general long-time existence criteria are still being developed \cite{PPZ18CAG}.

\subsection{Dimensional reductions} \label{dim-reduc}
The anomaly flow is defined over a possibly non-K\"ahler Calabi-Yau threefold. In this section, we describe two dimensional reductions of the anomaly flow: one to Riemann surfaces and another to $K3$ surfaces.

\begin{ex}
Reduction to a Riemann surface. The full details and results are from joint work with Z. Huang and T. Fei \cite{FHP2}. The geometric construction, due to \cite{Fei}, starts with a genus $g \geq 2$ Riemann surface $\Sigma$ together with a non-constant holomorphic map $\varphi: \Sigma \rightarrow \mathbb{P}^1$ satisfying $\varphi^* \mathcal{O}(2) = K_\Sigma$. We pullback a basis of $H^0(\mathbb{P}^1,\mathcal{O}(2))$ and denote the holomorphic $(1,0)$-forms on $\Sigma$ by $\mu_1$, $\mu_2$, $\mu_3$. The reference metric on the Riemann surface $\Sigma$ will be
\[
\omega_\Sigma = i \mu_1 \wedge \bar{\mu}_1 + i \mu_2 \wedge \bar{\mu}_2+i \mu_3 \wedge \bar{\mu}_3.
\]
We write $\varphi=(\alpha,\beta,\gamma)$ by writing the map using stereographic coordinates $\mathbb{P}^1 = S^2 \subseteq \mathbb{R}^3$. Let $(T^4, g_{T^4})$ be the 4-torus equipped with the flat metric, which we view as a hyperk\"ahler manifold and let $I,J,K$ be a triple of compatible complex structures satisfying $IJK=-id$. The associated K\"ahler forms are denoted $\omega_I$, $\omega_J$, $\omega_K$. We take the smooth 6-manifold to be the product $X = T^4 \times \Sigma$, however the complex structure is determined by the twistor construction
\[
J = j_\Sigma \oplus (\alpha I + \beta J + \gamma K).
\]
The complex manifold $X$ does not admit a K\"ahler structure. It does admit a holomorphic volume form $\Omega$, so it is a non-K\"ahler Calabi-Yau threefold. The ansatz
\[
\omega(t) = e^{2 f(t)} \omega_\Sigma + e^{f(t)} (\alpha \omega_I + \beta \omega_J + \gamma \omega_K), \quad f: \Sigma \rightarrow \mathbb{R}
\]
is conformally balanced and preserved by the anomaly flow. We take the gauge bundle to be a trivial bundle with $F=0$, so that the ${\rm Tr} \, F \wedge F$ term drops out from the flow equation. With this setup in place, the flow reduces to the following nonlinear PDE on a Riemann surface:
\begin{equation} \label{af-riem}
\partial_t e^f = (g_\Sigma)^{z \bar{z}} \partial_z \partial_{\bar{z}} u - \kappa u 
\end{equation}
where
\[
u = e^f + {\alpha' \over 2} \kappa e^{-f}.
\]
The function $\kappa$ is defined by $\kappa \omega_\Sigma = - \varphi^* \omega_{FS}$, hence $\kappa \leq 0$ with zeroes at the branch points of $\varphi$. The main result of \cite{FHP2} shows that the flow exists for all time $t \in [0,\infty)$ for initial data satisfying $u(x,0) \geq 0$. The flow does not converge however, as there are no stationary points in the regime $u \geq 0$, and $\lim_{t \rightarrow \infty} | \Omega |_{\omega(t)} \rightarrow 0$. After normalization by $\int_X |\Omega|_\omega
\, \omega^3$, the flow collapses the fibers and converges to a limiting metric on $\Sigma$.

On the other hand, for very small initial data, there is a finite-time blow-up where $\lim_{t \rightarrow T} | \Omega |_{\omega(t)} \rightarrow \infty$. It remains to extend the analysis of the PDE \eqref{af-riem} to a broader class of initial data, especially to the intermediate regime where stationary points are present. The solutions found in \cite{FHP} depend on the configuration of the zero points of the function $\kappa$; explicit solutions are found when the image under $\varphi$ of the zeroes of $\kappa$ all lie in a hemisphere. It is not known how the flow detects and reacts to the configuration of zeroes of $\kappa$.
\end{ex}

\begin{ex}
Reduction to a $K3$ surface. Details on the geometric construction can be found in \cite{FY}, and results on the flow are joint work with D.H. Phong and X.-W. Zhang \cite{PPZFY}. Let $S$ be a $K3$ surface with K\"ahler Ricci-flat metric $\omega_{K3}$. Let $X$ be the total space of a $U(1) \times U(1)$ principle bundle, with each $U(1)$ bundle equipped with a connection 1-form $\theta_i$ with $d \theta_i = \omega_i$ and $\omega_i \in 2 \pi H^2(S,\mathbb{Z})$ anti-self-dual with respect to $\omega_{K3}$. The manifold admits a complex structure under which $\theta = \theta_1 + i \theta_2$ is a $(1,0)$-form. The manifold also admits a holomorphic volume form $\Omega = \Omega_S \wedge \theta$, but cannot admit a K\"ahler structure. The ansatz
\[
\omega(t) = e^{u(t)} \omega_{K3} + i \theta \wedge \bar{\theta}, \quad u: S \rightarrow \mathbb{R}
\]
is conformally balanced and preserved by the anomaly flow. The gauge bundle is taken to be $\pi^* E_S \rightarrow X$ where $E_S \rightarrow S$ is a stable bundle of degree zero over $(S,\omega_{K3})$. The curvature $F$ is taken to be a fixed Hermitian-Yang-Mills curvature form $F \wedge \omega_{K_3}=0$ and does not flow. With this ansatz, the anomaly flow becomes
\[
\bigg[ {d \over dt} e^u \bigg] \, \omega_{K3}^2= \bigg[ \Delta_{\omega_{K3}} e^u \bigg] \, \omega_{K3}^2 - \alpha' i \partial \bar{\partial} (e^{-u} \rho)  + {\alpha' \over 2}  (i \partial \bar{\partial} u)^2 + \mu
\]
where $\rho_{\mu \bar{\nu}} = {i \over 2} (\omega_1 - i\omega_2)^\alpha{}_\mu (\omega_1+i\omega_2)_{\bar{\nu} \alpha}$ and
\[
\mu = {\alpha' \over 4} ( {\rm Tr} \, F \wedge F - {\rm Tr} \, R_{\omega_{K_3}} \wedge R_{\omega_{K_3}}) - \omega_1^2 - \omega_2^2.
\]
We assume the integrability constraint $\int_S \mu = 0$. For large enough initial data $u_0 = \log M$ with $M \gg 1$, the main result of \cite{PPZFY} states that the flow exists for all time. Furthermore, the flow converges in the limit $t \rightarrow \infty$ to a stationary point, recovering the Fu--Yau \cite{FY} solution to the Strominger system.  As the anomaly flow preserves the balanced class, the result of \cite{PPZFY} produces a solution to the Strominger system in the cohomology class
\[
[ | \Omega |_\omega \omega^2] = M [\omega_{K3}^2] + 2 [ \omega_{K3} \wedge i \theta \wedge \bar{\theta}].
\]
It is expected, though not yet proved, that there is finite-time blowup of the form $\lim_{t \rightarrow T} |\Omega|_{\omega(t)} \rightarrow \infty$ for very small initial data $0<M \ll 1$, in the same way as for the anomaly flow over Riemann surfaces.
\end{ex}

\subsection{Flowing by the Riemannian Ricci tensor}
We reviewed in \S \ref{sec:eom} how the Strominger system produces special solutions to the equations of motion of heterotic string theory. It was noticed by Ashmore--Minasian--Proto \cite{AMP} that the anomaly flow can also be understood as a parabolic version of the equations of motion. To see this, we first apply the variational formula \eqref{variation-omega} to the flow \eqref{AF} to derive the induced flow on the metric:
\[
\partial_t \omega = {1 \over 2|\Omega|_\omega} \Lambda_\omega \bigg( i \partial \bar{\partial} \omega - \alpha' ({\rm Tr} \, R \wedge R - {\rm Tr} \, F \wedge F) \bigg).
\]
Combining this with identity \eqref{eom-ddbar} and letting
\[
  e^{-2 \Phi} = |\Omega|_\omega, \quad H=i(\partial-\bar{\partial})\omega
\]
as before, we obtain
\begin{align} \label{af-ricci}
  \partial_t g_{\alpha \bar{\beta}} &= {e^{2 \Phi} \over 2} \bigg[ - R_{\alpha \bar{\beta}} - 2 \nabla_\alpha \nabla_{\bar{\beta}} \Phi + {1 \over 4} H_{\alpha mn} H_{\bar{\beta}}{}^{mn} \bigg] \nonumber\\
  &+ \alpha'  {e^{2\Phi} \over 2} i \Lambda_\omega ( {\rm Tr} \, R \wedge R - {\rm Tr} \, F \wedge F)_{\alpha \bar{\beta}},
\end{align}
where $R_{mn}$ is the Riemannian Ricci curvature of the metric $g$, summation over $m,n$ denotes real indices and $\nabla$ is the Levi-Civita connection. Therefore the anomaly flow creates a flow of the metric tensor by the Riemannian Ricci curvature plus extra terms. 
As the metric remains conformally balanced along the flow, the co-closedness condition \eqref{div-Ho}
\begin{equation} \label{d-dagger-H}
d^\dagger (e^{-2\Phi} H)\equiv 0
\end{equation}
holds at all times.

These evolution equations can be compared with the renormalization group flow equations \cite{OSW} at zeroth order in $\alpha'$. The equation for the flow of the metric tensor $g_{ij}$ matches up to the overall conformal factor of $e^{2 \Phi}$. We should add that the equations for $H$ do not match up with \cite{OSW}, as the equations there are based on Type II string theory where $dH=0$ (and no vector bundle $E$ is present), while here we work in heterotic string theory where $dH \neq 0$ due to anomaly cancellation and the vector bundle $E$ is included.

In summary, as explained in \cite{AMP}, the right-hand side of the flow is driven by the equations of motion, and the co-closed condition \eqref{d-dagger-H} is automatically implemented. This flow can also be compared with the Heterotic-Ricci flow proposed in \cite{MMShah}, where the metric evolves by the same equation, $dH \neq 0$, but the heterotic Bianchi identity is preserved rather than the co-closed condition.

The anomaly flow can also be compared to the pluriclosed flow \cite{GenRicciFlow, MGFJS, StreetsTian1}. There the dilaton is not present and $\alpha'=0$. The pluriclosed flow is defined by
\begin{equation} \label{pluriclosed}
\partial_t\omega = \partial \partial^\dagger \omega + \bar{\partial} \bar{\partial}^\dagger \omega + {i \over 2} \partial \bar{\partial} \log \det g, \quad \partial \bar{\partial} \omega = 0.
\end{equation}
The pluriclosed flow on the Riemannian metric tensor is derived in \cite{StreetsTian2} to be
\begin{equation} \label{gen-ricci}
\partial_t g_{ij} = - R_{ij} - {1 \over 2} (L_{\theta^\#} g)_{ij} + {1 \over 4} H_{i mn} H_j{}^{mn},
\end{equation}
where $\theta = - J d^\dagger \omega$. Let us now compare this expression with the anomaly flow with $\alpha'=0$.

It is well-known (e.g. \cite{MGF-Survey}) that the conformally balanced non-K\"ahler condition can be written as $\theta = - d \log |\Omega|_\omega$. Therefore conformally balanced metrics solve
\[
(L_{\theta^\#} g)_{ij} = \nabla_i \theta_j + \nabla_j \theta_i = - 2 \nabla_i \nabla_j \log |\Omega|_\omega,
\]
and so the anomaly flow \eqref{af-ricci} with $\alpha'=0$ can be written
\[
\partial_t g_{ij} = {e^{2 \Phi} \over 2} \bigg[ - R_{ij} - {1 \over 2} (L_{\theta^\#} g)_{ij} + {1 \over 4} H_{i mn} H_j{}^{mn} \bigg].
\]
Therefore the two different flows of non-K\"ahler complex geometry, namely the anomaly flow \eqref{AF} and the pluriclosed flow \eqref{pluriclosed}, after neglecting a conformal factor and $\alpha'$-terms, lead to the same equation for the evolution of the metric tensors. The big difference is that the condition $dH \equiv 0$ is preserved by the pluriclosed flow, while in the anomaly flow instead the coclosed condition $d^\dagger (e^{-2 \Phi} H) \equiv 0$ is the condition to be preserved. In this sense, there appears to be a duality between the two flows.

\subsection{Anomaly flow without corrections}
\subsubsection{Higher dimensions} We note that the anomaly flow can also be studied after setting the $\alpha'$-terms to zero, and that the equation can be generalized to dimension $n \geq 3$. The flow equation is then 
\begin{equation} \label{af-alpha0}
\partial_t (|\Omega |_{\omega(t)} \, \omega(t)^{n-1}) = i \partial \bar{\partial} \omega(t)^{n-2},
\end{equation}
with initial data satisfying
\[
d (|\Omega|_{\omega_0} \, \omega_0^{n-1})=0
\]
on a compact complex manifold of dimension $n$ with fixed holomorphic volume form $\Omega$. It is shown in \cite{PPZ19} that this equation has the short-time existence and uniqueness property, and that fixed points of this flow must be K\"ahler Ricci-flat metrics. Thus the flow produces a deformation path from a conformally balanced metric towards a torsion-free Calabi-Yau structure. It was proved in \cite{BedVez} that conformally balanced metrics nearby a torsion-free Calabi-Yau structure flow into and converge to a K\"ahler Ricci-flat limit.

\begin{rk}
  The setup is reminiscent of a common setup for geometric flows in $G_2$ geometry. The flow \eqref{af-alpha0} takes an $SU(n)$ structure $(\omega, |\Omega|_\omega^{-1} \Omega)$ with $\omega$ conformally co-closed to a torsion-free $SU(n)$ structure. Instead of arising from an $SU(n)$ structure, the geometry in the $G_2$ setting is determined by a 3-form $\varphi$. There are well-known $G_2$ flows which start with initial data satisfying a partial integrability condition (e.g. $\varphi$ closed or co-closed) and flow to stationary points which are torsion-free; for example there is the Laplacian flow \cite{Bryant} and the co-flow \cite{KMT} (see e.g. \cite{Lotay} for a survey and \cite{Kari} for general theory of $G_2$ flows). For connections between well-known flows in $G_2$ and complex geometry, we refer to \cite{AMP,Caleb}.
\end{rk}

The flow \eqref{af-alpha0} could then potentially detect whether the underlying complex manifold admits a K\"ahler structure. For a possible application, we note the conjecture in \cite{FinoVezz} which states that if a compact manifold $X$ admits both a pluriclosed metric ($\partial \bar{\partial} \omega_1 = 0$) and a conformally balanced metric ($d(|\Omega|_{\omega_2} \, \omega_2^{n-1})=0$), then $X$ admits a K\"ahler metric $\omega_3$. It would be interesting to see whether the presence of an auxiliary pluriclosed metric $\omega_1$ could be used as a barrier to prevent finite-time singularities along the flow.

As a consequence of the above discussion, the flow \eqref{af-alpha0} must diverge on a non-K\"ahler manifold. We note two simple instances of this phenomenon.

\begin{ex}
Consider the complex Lie group $SL(2,\mathbb{C})$ and let $X = SL(2,\mathbb{C}) / \Lambda$ be a compact quotient. The complex manifold $X$ does not admit a K\"ahler structure (e.g. \cite{Wang}). Let $\{ e_1, e_2, e_3 \}$ be a left-invariant basis of holomorphic vector fields on $SL(2,\mathbb{C})$ with structure constants $[e_i,e_j]= \epsilon_{ijk} e_k$. The ansatz \cite{FeiYau}
\[
\omega(t) = \rho(t) \sum_{a=1}^3 i e^a \wedge \bar{e}^a, \quad \Omega = e^1 \wedge e^2 \wedge e^3,
\]
where $\rho(t)$ is a constant, is conformally balanced and preserved by the flow. The flow develops a singularity at a finite time $T< \infty$ (see \cite{cetraro} for an exposition) where
\[
\lim_{t \rightarrow T} |\Omega|_{\omega(t)} \rightarrow 0.
\]
We note that the anomaly flow with $\alpha'$-corrections included has also been studied on examples with Lie group symmetry, and we refer to \cite{PPZ18Con,Pujia,PujiaUg}.
\end{ex}

\begin{ex}
Our second example is an Iwasawa-type manifold inspired by the Fu--Yau ansatz \cite{FY}. We let $X = \mathbb{C}^3 / \mathbb{Z}^3[i]$, where the Gaussian integers $a,b,c \in \mathbb{Z}[i]$ act by
\[
(x,y,z) \mapsto (x+a,y+c, z+ \bar{a} y + b).
\]
It is well-known that $X$ does not admit a K\"ahler structure (e.g. \cite{cetraro}). The following forms descend to $X$:
\[
\Omega = dz \wedge dx \wedge dy, \quad \theta = dz - \bar{x} dy, \quad \omega_{T^4} = i dx \wedge d \bar{x} + i dy \wedge d \bar{y}.
\]
We can view the geometry of $X$ as a fibration $\pi: X \rightarrow T^4$ with $\pi(x,y,z) = (x,y)$. We start the flow with the ansatz
\[
\omega(t) = e^{u(t)} \omega_{T^4} + i \theta \wedge \bar{\theta},
\]
where $u: T^4 \rightarrow \mathbb{R}$ is a function on the base $T^4$, so that $u=u(x,y)$. Direct calculation \cite{cetraro} shows that $\omega(t)$ is conformally balanced and the ansatz is preserved by the flow. In fact, the flow reduces to the following equation for the scalar function:
\[
\partial_t e^u = {1 \over 2} (\Delta_{\omega_{T^4}} e^u + 1).
\]
Standard PDE theory implies that the flow exists for all time $t \in [0,\infty)$. Since $|\Omega|_\omega = e^{-u}$, then 
\[
\lim_{t \rightarrow \infty} |\Omega|_{\omega(t)} \rightarrow 0,
\]
since $e^u \rightarrow \infty$. 
\end{ex}

Returning to the general theory, the estimates in \cite{FPPZ-B} show that singularities can only form if one of the estimates
\[
\omega(t) \leq C_1 \omega_0, \quad |\nabla \log |\Omega|_{\omega(t)} | \leq C_2
\]
fail. The upper bound on the metric implies a lower bound for $|\Omega|_\omega$. We wonder whether the a priori bound on $\nabla \log |\Omega|_\omega$ is needed, and whether the upper bound $\omega \leq C \omega_0$ can be replaced by the estimate $|\Omega|_{\omega(t)} \geq C^{-1} > 0$.

\subsubsection{Type IIA/IIB Flows} Finally, we note two flows related in spirit to equation \eqref{af-alpha0}. See \cite{Phong} for a more indepth survey of geometric flows from string theory.
\smallskip
\par $\bullet$ The first flow comes from adding a source term to the right-hand side. Let $X$ be a complex manifold of dimension $n \geq 3$ with holomorphic volume form $\Omega$. The Type IIB flow, introduced in \cite{FPPZ-B}, is defined by: 
\begin{align} \label{typeiib}
 & \partial_t (|\Omega|_\omega \omega^{n-1}) = i \partial \bar{\partial} \omega^{n-2} - \Psi, \nonumber\\
 & d (|\Omega|_{\omega_0} \, \omega_0^{n-1}) = 0
\end{align}
where $\Psi \in \Lambda^{n-1,n-1}(X,\mathbb{R})$ is a given closed form. Results for the Type IIB flow include \cite{FeiPhongB,FPPZ-B,Klemy}.
\smallskip
\par $\bullet$ There is a mirror to equation \eqref{typeiib} which is suggested by dualities in string theory, and there is a sense in which the dual of non-K\"ahler Calabi-Yau complex geometry is non-K\"ahler Calabi-Yau symplectic geometry. For further details on this, we refer to Tseng--Yau \cite{TsengYau}. To setup the associated geometric flow, let $(X,\omega)$ be a compact 6-dimensional symplectic manifold. The Type IIA flow, introduced in \cite{FPPZ-A}, is for a positive 3-form $\varphi(t)$:
\[
\partial_t \varphi = d \Lambda_\omega d( |\varphi|^2 \star \varphi) - \rho
\]
with initial data $\varphi_0$ a primitive, closed, positive 3-form \cite{Hitchin}. Here $\rho \in \Lambda^{3}(X,\mathbb{R})$ is a given closed form. Hitchin's construction \cite{Hitchin} produces from a positive 3-form $\varphi$ an almost-complex structure $J_\varphi$. Thus $(X,\omega, J_{\varphi(t)})$ is a triple where the symplectic form is fixed and $d \omega = 0$, but the almost-complex structure $J_{\varphi(t)}$ evolves and is non-integrable in general. For further developments on the Type IIA flow and applications to symplectic geometry, we refer to \cite{FeiPhong,FPPZ-A2,Klemy2,Raff}.

\appendix
\section{Conventions} \label{appendix}
Let $X$ be a complex manifold. We will let Roman letters $i,j,k$ denote real indices and Greek letters $\alpha,\beta,\gamma$ denote holomorphic indices, so that summation over $i$ includes summation over $\alpha$ and $\bar{\alpha}$. A hermitian form will be denoted
\[
\omega = i g_{\mu \bar{\nu}} \, dz^\mu \wedge d \bar{z}^\nu,
\]
where $g_{\mu \bar{\nu}}$ is an $n \times n$ positive-definite local hermitian matrix: $\overline{g_{\mu \bar{\nu}}} = g_{\nu \bar{\mu}}$. Its inverse is denoted $g^{\bar{\sigma} \mu}$ so that $g^{\bar{\sigma} \mu} g_{\mu \bar{\nu}} = \delta^{\bar{\sigma}}{}_{\bar{\nu}}$. Let $\nabla$ denote a covariant derivative on
\[
  T_{\mathbb{C}} X = T^{1,0}X \oplus T^{0,1} X,
\]
acting on vector fields
\[
  V = V^k \partial_k = V^\gamma \partial_\gamma + V^{\bar{\gamma}} \partial_{\bar{\gamma}}.
\]
Then $\nabla_i V = (\nabla_i V^k) \partial_k$ where
\[
\nabla_i V^k = \partial_i V^k + \Gamma_i{}^k{}_\ell V^\ell.
\]
The curvature tensor of $\nabla$ is denoted
\begin{equation} \label{curvature}
R_{pq}{}^m{}_j = \partial_p \Gamma_q{}^m{}_j + \Gamma_p{}^m{}_r \Gamma_q{}^r{}_j - (p \leftrightarrow q).
\end{equation}
With these conventions, when $\Gamma$ is the Levi-Civita connection then the Ricci tensor is given by
\begin{equation} \label{ricci}
R_{pj} = - R_{pm}{}^m{}_j.
\end{equation}
We now give specific examples of connections on $T_{\mathbb{C}} X$ used in the main text. The non-K\"ahler corrections will involve the following 3-form
\[
H = i (\partial-\bar{\partial}) \omega.
\]
Our conventions for components of a differential form $\eta \in \Lambda^{p,q}(X)$ are
 \[
\eta = {1 \over p! q!} \eta_{\alpha_1 \cdots \alpha_p \bar{\beta}_1 \cdots \bar{\beta}_q} dz^{\alpha_1} \wedge \cdots \wedge dz^{\alpha_p} \wedge d \bar{z}^{\beta_1} \wedge \cdots \wedge d \bar{z}^{\beta_q},
\]
so that for example
\[
H_{\alpha \beta \bar{\sigma}} = - (\partial_\alpha g_{\beta \bar{\sigma}} - \partial_\beta g_{\alpha \bar{\sigma}}).
\]
The Chern connection is defined by
\[
\Gamma^{\rm Ch}{}_\mu{}^\kappa{}_\nu = \partial_\mu g_{\nu \bar{\sigma}} g^{\bar{\sigma} \kappa},
\]
with
\[
\Gamma^{\rm Ch}{}_{\bar{\mu}}{}^{\bar{\kappa}}{}_{\bar{\nu}}  = \overline{\Gamma^{\rm Ch}{}_\mu{}^\kappa{}_\nu}
\]
and all other connection coefficients $\Gamma_i{}^k{}_j$ are zero. The conventions are such that
\[
  \Gamma^{\rm Ch}{}_\mu{}^\kappa{}_\nu - \Gamma^{\rm Ch}{}_\nu{}^\kappa{}_\mu = H_\mu{}^\kappa{}_\nu.
\]
The Levi-Civita connection is
\[
\Gamma_i{}^k{}_j = {g^{k \ell} \over 2}(-\partial_\ell g_{ij} + \partial_i g_{\ell j} + \partial_j g_{\ell i}),
\]
which in holomorphic coordinates is
\bea \label{LC2C}
\Gamma^{\rm LC}{}_\mu{}^\kappa{}_\nu &=& \Gamma^{\rm Ch}{}_\mu{}^\kappa{}_\nu - {1 \over 2} H_\mu{}^\kappa{}_\nu \nonumber\\
\Gamma^{\rm LC}{}_\mu{}^{\bar{\kappa}}{}_\nu &=& 0\nonumber\\
\Gamma^{\rm LC}{}_\mu{}^\kappa{}_{\bar{\nu}} &=& {1 \over 2} H_\mu{}^\kappa{}_{\bar{\nu}} \nonumber\\
\Gamma^{\rm LC}{}_{\mu}{}^{\bar{\kappa}}{}_{\bar{\nu}} &=& - {1 \over 2} H_{\mu}{}^{\bar{\kappa}}{}_{\bar{\nu}}. \nonumber
\eea
 The Hull connection in holomorphic coordinates is
\bea
\Gamma^{\rm H}{}_\mu{}^\kappa{}_\nu &=& \Gamma^{\rm Ch}{}_\mu{}^\kappa{}_\nu \nonumber\\
\Gamma^{\rm H}{}_\mu{}^{\bar{\kappa}}{}_\nu &=& 0\nonumber\\
\Gamma^{\rm H}{}_\mu{}^\kappa{}_{\bar{\nu}} &=& H_\mu{}^\nu{}_{\bar{\kappa}}\nonumber\\
\Gamma^{\rm H}{}_{\mu}{}^{\bar{\kappa}}{}_{\bar{\nu}} &=& 0.\nonumber
\eea
The Strominger-Bismut \cite{Bismut,Strominger} connection $\hat{\nabla}$ is defined on a section $V \in \Gamma(T_{\mathbb{C}} X)$ by
\begin{equation} \label{Str-Bis}
\hat{\nabla}_\mu V^\kappa = \nabla^{\rm Ch}_\mu V^\kappa + H^\kappa{}_{\mu \lambda} V^\lambda, \quad \hat{\nabla}_{\bar{\mu}} V^\kappa = \partial_{\bar{\mu}} V^\kappa - V^\lambda H_{\lambda \bar{\mu}}{}^\kappa
\end{equation}
and $\nabla_i V^{\bar{\kappa}} = \overline{\nabla_i V^\kappa}$. The connection $\hat{\nabla}$ was noticed by Yano \cite{Yano} as the only connection preserving the metric, the complex structure and with skew-symmetric torsion tensor.

Finally, in the main text we will use the expression for the curvature of the Hull connection. Recall our conventions are such that $R^{\rm H}{}_{pq}{}^m{}_n $ has 2-form indices $p,q$ and endomorphism indices $m,n$. The $(1,1)$ part of the Hull curvature is then
\bea \label{RH11}
R^{\rm H}{}_{\mu \bar{\nu}}{}^\alpha{}_\beta &=& R^{\rm C}{}_{\mu \bar{\nu}}{}^\alpha{}_\beta + H_\mu{}^\alpha{}_{\bar{\sigma}} H_{\bar{\nu}}{}^{\bar{\sigma}}{}_\beta \\
R^{\rm H}{}_{\mu \bar{\nu}}{}^{\bar{\alpha}}{}_\beta &=& \nabla_\mu H_{\bar{\nu}}{}^{\bar{\alpha}}{}_\beta \nonumber\\
R^{\rm H}{}_{\mu \bar{\nu}}{}^\alpha{}_{\bar{\beta}} &=& - \nabla_{\bar{\nu}} H_{\mu}{}^\alpha{}_{\bar{\beta}}. \nonumber
\eea
The $(2,0)$ part of the curvature works out to be
\bea \label{RH20}
R^{\rm H}{}_{\mu \nu}{}^\alpha{}_\beta &=& 0 \\
R^{\rm H}{}_{\mu \nu}{}^{\bar{\alpha}}{}_\beta &=& 0 \nonumber\\
R^{\rm H}{}_{\mu \nu}{}^\alpha{}_{\bar{\beta}} &=& - (i \partial \bar{\partial} \omega)_{\mu \nu}{}^\alpha{}_{\bar{\beta}}. \nonumber
\eea

\end{document}